\RequirePackage{fix-cm}
\documentclass{svjour3}                     
\smartqed  

\usepackage{graphicx}

\usepackage{amssymb}
\usepackage{amsmath}
\usepackage{mathrsfs} 
\usepackage[T1]{fontenc}

\usepackage{theorem}
\usepackage{fancyhdr}
\usepackage[mathscr]{euscript}
\usepackage[sans]{dsfont}
\usepackage{tikz}
\usetikzlibrary{shapes,arrows}
\usetikzlibrary{matrix}
\usepackage{graphicx}
\usepackage{color}

\usepackage{hyperref}

\newcommand{\e} {\varepsilon}
\newcommand{\C} {\mathbb{C}}
\newcommand{\D} {\mathbb{D}}

\newcommand{\Z} {\mathbb{Z}}
\newcommand{\N} {\mathbb{N}}

\newcommand{\R} {\mathbb{R}}

\newcommand{\inte}{ \operatorname{int}}

\newcommand{\di}[1]{\operatorname{d}\!#1}

\newcommand{\p} {\mathbf{P}}

\begin{document}

\title{Polynomial methods to  construct inputs for uniformly ensemble reachable  linear systems
}

\titlerunning{Constructive methods for ensemble reachability}        

\author{Michael Sch\"onlein
}


\institute{  Michael Sch\"onlein \at
              Bauhaus University Weimar \\
              Coudraystra\ss{}e 13\\
              99423 Weimar\\
              Germany\\
              \email{michael.schoenlein@uni-weimar.de}           
}

\date{Received: date / Accepted: date}

\maketitle

\begin{abstract}
This paper is concerned with linear parameter-dependent systems and considers the notion of uniform ensemble reachability. The focus of this work is on constructive methods to compute suitable parameter-independent open-loop inputs for such systems. In contrast to necessary and sufficient conditions for ensemble reachability,  computational methods have to distinguish between continuous-time and discrete-time systems. Based on recently derived sufficient conditions and techniques from  polynomial approximation, we present two methods for discrete-time {single}-input linear systems. Moreover, we illustrate that one method can also be applied to certain continuous-time single-input systems.

\end{abstract}

{}

\section{Introduction}\label{intro}

Ensemble control is an emerging field in mathematical systems and control theory  referring to the task of controlling a large, potentially infinite, number of states, or systems, using a single input function or a single feedback controller. Being not a {rigorously} defined term it subsumes several different scenarios, where in each particular case the problems and the required techniques differ from each other. Ensemble control embraces the scenario of uncertainties in the initial data. This is modeled by  a probability distribution of the initial state and the ensemble control problem leads to a transport problem of density functions and therefore to controllability and observability
issues of the Liouville and the Fokker-Plank equation
\cite{brockett2012notes}, \cite{chen2017optimal}, \cite{fleig2016estimates}, \cite{shen2017discrete}, \cite{zeng2016ensemble}.
Another setting tackles uncertainties in the model parameters by addressing controllability problems for parameter-dependent {systems,}
and the goal is to achieve a certain control task by using only a single or a few open-loop inputs which are \textit{independent} of the (usually unknown) model parameters. The notion \emph{ensemble controllability} is used for this, cf. \cite{li2009ensemble}. In this context we also mention the notions \emph{averaged controllability} and  \emph{moment controllability}, cf. \cite{zuazua2014averaged} and \cite{Zeng_scl_2016}, respectively. Motivated by the controllability analysis of the Bloch equation, there is also the concept of asymptotic controllability, cf. e.g. \cite{beauchard2010controllability}, \cite{Gauthier_scl_2018}.

In this paper we consider families of parameter-dependent systems  of the form
\begin{align} \label{eq:def-ensemble}\begin{split}
\tfrac { \partial } { \partial t } x ( t ,  \theta  ) &= A (  \theta ) x ( t ,  \theta  ) + B (  \theta  ) u ( t ) 
\\
x _{ t + 1 }( \theta ) &= A ( \theta ) x _ t ( \theta ) + B ( \theta ) u _t, \end{split}
\end{align}
where $\theta\in \p$ is considered as a parameter and the parameter space $\p\subset \C$ is assumed to be a Jordan arc in the complex plane, i.e. $\p$ is the image of a continuous and bijective function defined on a compact interval. Moreover, we assume that the matrix-valued functions $A: \p \to \C^{n \times n}$ and $B: \p \to \C^{n \times m}$ are continuous. Throughout the paper we will use the short notation  $(A,B) \in C_{n,n}(\mathbf{ P })\times C_{n,m}(\mathbf{ P })$. In addition, let $C_n(\p)$  denote the space of continuous functions from $\p$ to $\C^n$. As we mainly consider discrete-time systems, we take the initial {condition} $x(0,\theta)=x_0(\theta)=0$. To express that the solutions to \eqref{eq:def-ensemble} are regarded as functions from the parameter space $\p$ to $\C^n${,} we denote them by $\varphi(T,u)(\theta):=\varphi(T,u,0,\theta)$. Further, we wish to highlight the essential property that the input is {\it independent} of the parameters.

The notion of reachability we are considering in this paper is as follows: We say that a pair $(A,B) \in C_{n,n}(\mathbf{ P })\times C_{n,m}(\mathbf{ P })$  is \textit{uniformly ensemble reachable} (from zero), if for every $f \in C_n(\p)$ and $\e>0$ there are $T>0$ and  $u\in L^1([0,T],\C^m)$ or $u = (u_{0},u_{1}, \dots, u_{T-1}), u_{i} \in \C^{m}$ such that   $$ \|\varphi(T,u)-f\|_{\infty}:= \sup_{\theta \in \p} \|\varphi(T,u,0)(\theta)-f(\theta) \|_{\C^n} <\e.$$ 

Note that ensemble reachability is an infinite-dimensional problem (unless $\p$ is finite) and it is equivalent to approximate reachability of the infinite dimensional linear system defined by the bounded linear matrix-multiplication operators 
\[ {\cal M}_A\colon C_n(\p) \to C_n(\p), \quad {\cal M}_A f(\theta)=A(\theta)f(\theta)\]
and 
\[
 {\cal M}_B\colon \C^m \to C_n(\p),\quad {\cal M}_B u =B(\theta)u.\]
As a consequence of the restriction that the inputs are not allowed to {depend} on the parameter, it is well-known that exact reachability (i.e. $ \e=0$) is never possible provided $\p$ is not finite, cf. \cite[Theorem~3.1.1]{triggiani75} and \cite[p.244]{fuhrmann_hilbert_book}. {We note that in the literature
continuous-time systems are mostly considered, and hence
the term {\it ensemble controllability} is present more frequently.}
Moreover, it follows from \cite[Theorem~3.1.1, Remark~3.1.2]{triggiani75} that {a continuous-time pair} $(A,B)$ is uniformly ensemble reachable if and only if it is {\it completely approximately controllable}, i.e.{,} for every $T>0$, for every $\e>0$ and for every pair $x_0,x_1 \in C_n(\p)$ there exists {an input} $u \in L^1([0,T],\C^m)$ such that
\begin{align*}
 \| x_1 - \varphi(T,u,x_0)\|_{\infty} < \e. 
\end{align*}
{We} note that this equivalence does not hold for discrete-time systems.
Recall that a pair $(A,B) \in C_{n,n}(\p)\times C_{n,m}(\p)$  is called {\it uniformly ensemble controllable} (to zero), if for all $x_0 \in C_n(\p)$ and $\e >0$ there {exist} $T \geq 0$
and an input $u$ such that
\begin{align*}
\|\varphi(T,u,x_0)\|_{\infty} < \e.
\end{align*}
It is well-known that the notions of approximate reachability (from zero) and approximate controllability (to zero)
are independent of each other, cf.~\cite[Lemma~4.1]{fuhrmann1972weak}{,} and none of both implies approximate
complete controllability. {Also we note that this paper does not cover other notions of controllability, such as asymptotic controllabilty, cf. \cite{sontag_sicon_1983_asym_control}.} Moreover, recall that the application of \cite[Theorem~3.1.1]{triggiani75} yields that a pair $(A,B)$ is uniformly ensemble reachable if and only if
\begin{align}\label{eq:trig-control}
 \overline{\operatorname{span} \{\theta \mapsto A(\theta)^lb_j(\theta) \,\, |\, \, j=1,...,m, \, l=0,1,...\} } = C_n(\p),
\end{align}
where $b_j(\theta)$ denotes the $j$th column of $B(\theta)$. In operator theory language, uniform ensemble reachability of $(A,B)$ is  equivalent to the fact that the matrix multiplication operator ${\cal M}_A$ is $m$-multicyclic and the columns $b_1,...,b_m$ are cyclic vectors. For other characterizations{,} we refer to \cite[Section~6.2]{curtainzwart95} and \cite[Section~VI.8.a]{engelnagel}. As these characterizations have the drawback that they are hard to check in practice, much effort has recently been spent on the derivation of pointwise testable necessary and sufficient conditions, cf. \cite{JDE_ensembles_2021}, \cite{li_tac_2016}.

The ensemble controllability problem for continuous-time parameter-dependent systems is also studied in \cite{Agrachev_ensemble_lie_2016}, \cite{agrachev2020control}, \cite{chen2019mcss}, \cite{li2011}, \cite{li2009ensemble}, \cite{li_tac_2016}, \cite{loheac:hal-01164525} and  \cite{Zeng_scl_2016}. Agrachev and Sarychev consider ensemble controllability for nonlinear drift-free parameter-dependent systems and provide a characterization in terms of Lie-brackets. In the same direction, the work of Chen \cite{chen2019mcss} also treats nonlinear systems and considers Lie extensions. We note that these approaches do not apply to the setting in this paper. In \cite{li2011} a characterization for ensemble controllability for time-varying parameter-dependent linear systems is presented, which is based on the singular value decomposition of the reachability operator. Like condition \eqref{eq:trig-control}, this condition is hard to check. The recent works \cite{chen_multi_arxiv} and \cite{danhane:hal-03824645} have shown that, if $(A,B)$ satisfy appropriate regularity assumptions, $L^p${-} and uniform ensemble controllability, respectively, cannot hold if the parameter space $\p\subset \R^d$, $d \geq 2$ has interior points.

\medskip

{\it Problem statement} $ \,$ The objective of this paper is as follows. We consider uniformly ensemble reachable pairs $(A,B) \in C_{n,n}(\p) \times C_{n,m}(\p)$. Given a desired target family $f \in C_n(\p)$ and an $\e$-neighborhood of it, we treat the problem {of} how to find a suitable $T>0$ and a suitable input $u$ such that
\begin{align*}
 \|\varphi(T,u) - f\|_{\infty} < \e.
\end{align*}

For continuous-time systems this problem has recently been addressed in \cite{lazar:hal-03035494}, \cite{scindis} {and \cite{Zeng_scl_2022}}. {
The methods in \cite{scindis} and \cite{Zeng_scl_2022} are based on techniques from integral equations}. Also we mention that the results in \cite{scindis} require an additional assumption on the target family $f$ so that predefined error bounds can be archived. Besides {that,} the work \cite{lazar:hal-03035494} uses the adjoint system to set up a continuous-time optimal control problem to derive suitable inputs. However, as the solution formula for discrete-time systems is not given in terms of an integral, the methods just mentioned cannot be applied to discrete-time systems. Moreover, we are not aware of any approach{es} in the literature to tackle the discrete-time {case}.

\medskip
\noindent
{\it Novelty and main contribution} $\,$
The scope of the present paper is to fill this gap. {More precisely} we focus on single-input systems
and start with treating the discrete-time case 
\begin{align*}
x _{ t + 1 }( \theta ) = A ( \theta ) x _ t ( \theta ) + {b} ( \theta ) u _t.
\end{align*}
The solution at time $T>0$ for inputs $u=(u_0,...,u_{T-1}) \in \C \times \cdots \times \C$ is given by
\begin{equation*}
\varphi(T,u)(\theta) =  
\sum_{k=0}^{T-1} \big(A(\theta)\big)^{T-1-k}b(\theta)u_{k}.
\end{equation*}
Noticing the restriction that the inputs $u_0,...,u_{T-1}$  have to be independent of the parameters, the solution  can be expressed in terms of the polynomial $$ p : \C \to \C,\quad p(z)= u_{T-1} + u_{T-2}\,z+\cdots + u_{1}\, z^{T-2}+ u_{0} \,z^{T-1}.$$
That is, we have
 \[
 \varphi(T,u)(\theta) = p(A(\theta))b(\theta).
 \]
In this notation, the problem can be stated equivalently as follows. Given $f\in C_n(\p)$ and $\e>0${,} how to find a suitable polynomial $p$ such that
 \begin{align}\label{eq:key-obs}
  \| \varphi(T,u) - f\|_{\infty} = \|p(A)b -f\|_\infty < \varepsilon \, {?}
 \end{align}
This observation is the key to a new strategy to solve the construction problem by using tools from polynomial approximation. In this context, given the target family $f \in C_n(\p)$ and the $\e$-neighborhood of it, the task of finding the amount required inputs and their values is equivalent to derive the degree and the coefficients of an approximating polynomial (in the sense of \eqref{eq:key-obs}). Note that this fact is significantly different from the continuous-time case, where the solution formula {is given in terms of} an integral operator.\\
\indent
The contributions of this paper are made up of two parts. On {the} one hand, the paper provides a collection of results from (complex) polynomial approximation, for which explicit representations of the  corresponding polynomials are available. More precisely, after a short recap of explicit refinements of the classical Weierstrass Approximation Theorems,  the paper provides an explicit construction of the polynomial in Runge's Little Theorem, where the degree and the coefficients are determined in terms of the target function and the approximation error.  To the best of the author's knowledge, such an explicit representation was not given in the literature so far. Moreover, the paper gives an explicit treatment of results which are due {to} Walsh. As these results are required later for the derivation of suitable inputs, the paper also provides a detailed treatment of these results (which is only implicitly contained in the literature). More specifically, we show that the construction of {appropriate} polynomial{s} can be traced back to Runge's Little Theorem by using suitable conformal mappings. We note that the numerical treatment of conformal mappings {is} a wide area. Hence, this topic is beyond the scope of this paper and we only provide some comments to {the} relevant literature. \\
\indent
On {the} other hand, based on the exposition on polynomial approximation, the paper shows various methods to determine appropriate inputs. As there is no explicit characterization of uniform ensemble reachability\footnote{The characterization in \eqref{eq:trig-control}{,} is not explicit in terms of $(A,B)$.} we take two sufficiency conditions, denoted by $(S1)$ and $(S2)$ {in Section~3}, derived in \cite{JDE_ensembles_2021} as a starting point. In this paper we  present constructive methods for each condition. Moreover, the paper shows that the  methods corresponding to the sufficiency condition $(S2)$ can also be applied to continuous-time single-input systems, while the methods for $(S1)$ are limited to discrete-time systems.

 \medskip

{\it Organization} $ \,$ 
{Let's  quickly describe the structure of the paper}.  Section~\ref{sec:approx} contains the relevant results from polynomial approximation that are required for the input construction procedures. It is divided into three subsections corresponding to the results of Runge, Walsh and Weierstrass, respectively. Section~\ref{sec:poly_si_methods} is devoted to the presentation of the constructive procedures to determine suitable inputs for discrete-time parameter-dependent linear single-input systems. After a short summary of relevant known results, the section {has} two subsections, where each subsection presents methods corresponding to one of the sufficiency conditions. Section~\ref{sec:poly_approx_cont} depicts that methods corresponding to sufficiency condition $(S2)$ can also applied to continuous-time systems. Finally, {Section~\ref{sec:comments}} contains additional perspectives about the presented methods.

 \medskip

{
{\it Notations} $ \,$
For $x\in \R$, let $\lfloor x \rfloor:= \max\{m \in \Z \, :\, m \leq x\}$.
For a set $\Omega$, its interior is denoted by $\inte \Omega$ and its closure is denoted by $\overline{\Omega}$.   The distance between two sets $\Omega,\tilde\Omega \neq \emptyset$ is denoted by
$$ \di{\,(\Omega,\tilde \Omega)} := \inf\{ \|\omega - \tilde \omega\| \, : \, \omega \in \Omega, \, \tilde \omega \in \tilde\Omega\}.$$
Moreover, we call $C_n(\Omega)$ the set of continuous functions $g \colon \Omega \to \C^n$ and we say that $g \in C_n^k(\Omega)$ if $g$ is $k$-times continuously differentiable. A function $g \colon \Omega \to \C^n$ is said to satisfy a Lipschitz condition if there is a $L_g>0$ such that
\begin{equation*}
 \|g(z_1) - g(z_2)\|_{\C^n} \leq L_g \, \|z_1- z_2\|_{\Omega}
\end{equation*}
for all $z_1, z_2 \in \Omega$. The set of functions satisfying a Lipschitz condition is denoted by $\operatorname{Lip}_n(\Omega)$. If $\Omega \subset \C$ is compact and  $g \colon \Omega \to \C$ is continuous, we define  $$M_{g} := \max_{z \in \Omega} |g(z)|.$$
For a function $f\colon [a,b] \to \R$, the $n$th Bernstein polynomial is given by
\begin{align}\label{def:poly_Bernstein}
 B_{n,f}(x):=   \frac{1}{(b-a)^n}\sum_{k=0}^n {n \choose k}\, f(a+\tfrac{k}{n}(b-a)) \,(x-a)^k  \, (b-x)^{n-k}.
\end{align}
For a continuous function $f\colon \partial \mathbb{D} \to \C$, its $k$th \emph{Fourier coefficient} is denoted by
\begin{align}\label{def:poly_Fejer_coeff}
\hat f(k) :=  \frac{1}{2 \pi} \int_{-\pi}^\pi f(e^{is}) e^{-iks} \di{s}.
\end{align}
Further, we call
\begin{align}\label{def:poly_Fejer}
F_{f,n}(z) := \sum_{k=-n+1}^{n-1} \frac{n-|k|}{n} \hat f (k)\, \,  z^k,
\end{align}
the $n$th Fej\'{e}r polynomial.
}
\section{Elements of Approximation Theory}\label{sec:approx}

In this section we discuss results from approximation theory{,} for which constructive proofs are available. The presented results  are due to Bernstein, Runge, Weierstrass and Walsh. Note that the definite result in this context is the famous Mergelyan's Theorem{,} saying that a continuous function $f \colon K \to \C$ can be uniformly approximated by polynomials if $K$ is compact, $\C\setminus K$ is connected and $f$ is analytic in  the interior of $K$, cf. \cite[Chap.~III, \S~2, Theorem~1]{gaier1987}. However, its highly ingenious method of proof is not constructive. Therefore we will present some special cases which can be {proved} constructively and prepare the ground for the computational methods in Section~\ref{sec:poly_si_methods}. {In} Section~\ref{sec:approx_weier} we recall known explicit versions of the classical Weierstrass approximation theorems. Moreover, in Section~\ref{sec:approx_runge} we present an explicit treatment of the approximating polynomial in Runge's Little Theorem. In Section~\ref{sec:appox_walsh} {we} present a detailed derivation of some of Walsh's results on complex approximation. Note that these results appeared more then twenty years before Mergelyan's result.

\subsection{Weierstrass Theorems}\label{sec:approx_weier}

We start with the Weierstrass approximation theorems and consider the special cases where the continuous function that shall be approximated additionally satisfies a Lipschitz condition. Recall that for a function $f\in C([a,b],\R)$ the $n$th Bernstein polynomial is given by  
\begin{align*}
 B_{n,f}(x):=   \frac{1}{(b-a)^n}\sum_{k=0}^n {n \choose k}\, f(a+\tfrac{k}{n}(b-a)) \,(x-a)^k  \, (b-x)^{n-k}.
\end{align*}
{For} complex-valued functions $f \colon [a,b] \to \C$ one can apply the latter to the real and imaginary part. That is, for $f(x) = g(x)+ih(x)$ we consider the complex Bernstein polynomial $B_{n,f}(x) = B_{n,g}(x) + i B_{n,h}(x) $. Next, we recap a result of Gzyl and Palacios that provides an explicit error bound, cf. \cite{gzyl1997}.

\begin{theorem}[Gzyl, Palacios (1997)]\label{thm:gzyl}
Let $f\colon [a,b]\to \C$ satisfy a Lipschitz condition. Then, for $n \geq 3${,} the sequence of complex Bernstein polynomials satisfies
\begin{align*}
\|f -  B_{n,f}\|_{\infty} \leq \sqrt{2}
\left( 4\,M_f + \frac{(b-a)\,L_f}{2} \right) \frac{\sqrt{\log n}}{\sqrt{n}}.
\end{align*}
\end{theorem}
Note that, since Gyzl and Palacios consider real-valued functions defined on the unit interval, the constants in the previous statement are adjusted. Similarly, the second or trigonometric Weierstrass approximation theorem can also be proven constructively. It is based on Fej\'{e}r's Theorem and we first {recall their definition. Let $f\colon \partial \mathbb{D} \to \C$ be a continuous function, then the the $n$th Fej\'{e}r polynomial is given by}
\begin{align*}
F_{f,n}(z) := \sum_{k=-n+1}^{n-1} \frac{n-|k|}{n} \hat f (k)\, \,  z^k,
\end{align*}
where {$\hat f(k) $}
denotes the $k$th Fourier coefficient of $f$. From \cite[Ch.~VIII, Sec.~1, Thm.~1 and Sec.~2, Thm.~3]{natanson1964constructive} we recap the following result.

\begin{theorem}[Second Weierstrass Theorem]\label{thm:weier-circ}
Suppose that $f\colon \partial \mathbb{D} \to \C$ is continuous. Let  $(F_n)_{n\in \N}$  denote the  Fej\'{e}r polynomials.
\begin{enumerate}
 \item[(a)]  {The} sequence $(F_n)_{n\in \N}$ converges uniformly to $f$.
 \item[(b)] If $f$ satisfies a Lipschitz condition{,} it holds
\begin{align*}
 \sup_{z \in \partial \mathbb{D}} |f(z)- F_{f,n}(z)| \leq 2 \sqrt{2} \pi \,L_f  \cdot \frac{\ln n }{n},  
\end{align*}
where $L_f>0$ denotes the Lipschitz constant.
\end{enumerate}
\end{theorem}

\medskip
We note that sharper but less explicit error bounds for the second Weierstrass approximation Theorem have been derived in L.~Lorch~\cite{lorch1962approximation}, S.M.~Nikolski~\cite{nikolski1940allure} and S.A.~Telyakovskii \cite{telyakovskii1969approximation}.

\subsection{Runge's Little Theorem}\label{sec:approx_runge}

Another famous result from approximation theory is due to Runge. Before presenting a constructive proof, we have to fix some notation. The presentation of Runge's little Theorem is based on \cite{greene_krantz_1997} and \cite{remmert2013classical}. Let $\gamma$ denote a closed (piecewise) $C^1$-path. With a slight abuse of notation{,} we denote its trace also by  $\gamma$. Then,
\begin{align*}
 \operatorname{ind}_{\gamma}(z) := \frac{1}{2\pi i} \int_{\gamma} \frac{1}{\xi-z} \di{\xi}\,,
\quad z\in \C\setminus \gamma
\end{align*}
denotes the winding number of $z$ with respect to $\gamma$. Moreover, a closed polygon $\tau=[p_1\,p_2\,\cdots\, p_k\, p_1]$ 
composed of finitely many horizontal or vertical segments $[p_1\,p_2]$,$[p_2\,p_3]$,...,$[p_k\,p_1]$ is called a \emph{grid polygon} if
there exists a not necessarily regular grid $G \subset \C$ of horizontal or vertical lines such that all vertices $p_1,...,p_k$ are pairwise distinct adjacent grid {points} of $G$. Moreover, let
\begin{align*}
 \operatorname{ext} \gamma := \{ z\in \C\setminus \gamma \, | \,  \operatorname{ind}_{\gamma}(z) =  0\}
\end{align*}
and
\begin{align*}
 \operatorname{int} \gamma := \{ z\in \C\setminus\gamma\, | \,  \operatorname{ind}_{\gamma}(z) =  1\}
\quad \text{ or } \quad  
\operatorname{int} \gamma := \{ z\in \C\setminus\gamma\, | \,  \operatorname{ind}_{\gamma}(z) =  -1\}
\end{align*}
depending on the orientation of $\gamma$. We will not provide a complete proof. Instead, as we are interested {in} constructive methods, we will provide an explicit treatment of the degree of the approximating polynomial as well as an explicit representation of its coefficients.  To best of the {author's} knowledge these are not available in the literature. For some analytical parts that will be left out we refer to \cite{remmert2013classical}.

\begin{theorem}[Runge's Little Theorem (1885)]\label{thm:runge}
Let $K\subset \C$ be compact such that $\C \setminus K$ is connected. If there is an open set $\Omega$ containing $K$ such that $f$ is holomorphic on $\Omega$, then for every $\e>0$ there is {a} polynomial $p$ such that
\begin{align*}
 \sup_{z\in K} |f(z)-p(z)|< \e.
\end{align*}
\end{theorem}
\begin{proof}
The proof is carried out in four constructive steps. We start with the construction of finitely many grid polygons in $\Omega\setminus K$ so that $f$ can be represented via the Cauchy integral formula for compact sets. 

\medskip
\noindent
{\bfseries Step~1: Grid polygon construction} \emph{There are finitely many distinct oriented horizontal or vertical line segments $\tau_1,...,\tau_N$ of length $\delta < \tfrac{1}{\sqrt{2}} \operatorname{d}(K,\partial \Omega)$ so that 
\begin{align*}
 f(z)= \sum_{k=1}^N    \frac{1}{2\pi i} \int_{\tau_k} \frac{f(\xi)}{\xi-z} \di{\xi}\qquad \forall \, \, z \in K.
\end{align*}
}
\medskip
\noindent
\emph{Proof of Step~1:}
Let  $\delta\in (0,\tfrac{1}{\sqrt{2}} \operatorname{d}(K,\partial \Omega))$ and consider a grid consisting of lines that are parallel to the real and imaginary axis so that the distance between two parallel lines is $\delta$. Since $K$ is compact{,} there are finitely many boxes $Q_1,...,Q_n$ that intersect with $K$ and satisfy
\begin{align*}
 K \subset \bigcup_{k=1}^n Q_k \subset \Omega.
\end{align*}
Indeed, for $q_l \in Q_l$ it follows from $2\delta < \operatorname{d}(K,\partial \Omega)$ that $B_{\delta}(q_l) \subset \Omega$. As each box has diameter $\sqrt{2}\delta$ and 
\begin{align*}
\sqrt{2}\delta < \operatorname{d}(K,\partial \Omega)<  \operatorname{d}(q_l,\partial \Omega){,}
\end{align*}
it follows that $Q_l\subset \Omega$.
 
The boundaries of the boxes $Q_1,...,Q_n$ consist of line segments. We now pick those boundary segments that are not common for two distinct boxes $Q_l$ and $Q_m$ with $l\neq m$. Let $\tau_1,...,\tau_N$ denote the selected segments. The construction yields that
\begin{align*}
 \bigcup_{k=1}^N \gamma_k \subset \Omega\setminus K.
\end{align*}
Indeed, if line segment $\tau_k$ meets $K$, there are two boxes of the grid having $\tau_k$ as a common side, a contradiction to the choice of the segments $\tau_1, ..., \tau_N$.

Let $z \in K$ be a point that does not lie on any boundary of the boxes $Q_1,...,Q_n$. {Then} there is exactly one box $Q_l$ containing $z$ and it holds
\begin{align*}
 f(z)=  \frac{1}{2\pi i} \int_{\partial Q_l} \frac{f(\xi)}{\xi-z} \di{\xi} \quad \text{ and } \int_{\partial Q_k} \frac{f(\xi)}{\xi-z} \di{\xi} =0 \quad \forall\,  k\neq l .
\end{align*}
Since common sides of different boxes of the grid occur twice with different orientation{,} it follows that
\begin{align*}
 f(z)=  \frac{1}{2\pi i}  \sum_{k=1}^n \int_{\partial Q_k} \frac{f(\xi)}{\xi-z} \di{\xi} = \frac{1}{2\pi i}  \sum_{k=1}^N \int_{\gamma_k} \frac{f(\xi)}{\xi-z} \di{\xi} =0.
\end{align*}
Moreover, as shown in \cite[Chapter~12, \S~1.1]{remmert2013classical}, the last equality also holds for all $z \in K$ that are on the boundary of some box.
$\,$ \hfill {\small $\Box$}

\medskip
Now, let $\e>0$, $\delta>0$ and $N$ be fixed. Moreover, let $L>0$ denote the Lipschitz constant of the function $\xi \mapsto  \tfrac{f(\xi)}{\xi-z}$ on $ \tau$.

\medskip

\noindent
{\bfseries Step~2:  Rational approximation} \emph{Let $M:= \lfloor\tfrac{3}{\e}\tfrac{2\pi}{\delta^2 \,L\,N}\rfloor$ and divide each line segment in $M$ subsegments of length $\tfrac{\delta}{M}$.  Let  the subsegments be denoted by $\tau_{k,l}$, for $l=1,...,M$ and $k=1,...,N$. Take distinct points  $w^{(k)}_{l} \in  \tau_{k,l}$ such that $|w^{(k)}_{l+1} - w^{(k)}_{l}| \leq \tfrac{\delta}{M}$ and let
\begin{align*}
c^{(k)}_l:= f(w^{(k)}_l) \int_{\tau_{k,l}}  \di{\xi} \quad \text{ for }  l=1,...,M \text{ and } k=1,...,N.
\end{align*}
Then, the rational function
\begin{align*}
r(z)= \frac{1}{2\pi i} \sum_{k=1}^{N} \sum_{l=1}^{M}  \frac{c^{(k)}_l}{w^{(k)}_l-z}
\end{align*}
approximates $f$ uniformly on $K$, i.e.{,} $r$ satisfies
\begin{align*}
 \sup_{z\in K} |  f(z)-r(z)|< \tfrac{\e}{3}.
\end{align*}
}
\medskip
\noindent
\emph{Proof of Step~2:} 
It holds
\begin{align*}
\left| \int_{\tau_k} \frac{f(\xi)}{\xi-z}\di{\xi}  -  \sum_{l=1}^{M}  \frac{c^{(k)}_l}{w^{(k)}_l-z} \right| 
 &=
 \left| \sum_{l=1}^M \int_{\tau_{k,l}} \frac{f(\xi)}{\xi-z}\di{\xi}  
 -  \sum_{l=1}^{M} \int_{\tau_{k,l}} \frac{f(w^{(k)}_l)}{w^{(k)}_l-z} \di{\xi}
 \right| 
\\
&\leq \sum_{l=1}^M
\left|  \int_{\tau_{k,l}} \left( \frac{f(\xi)}{\xi-z}  -    \frac{f(w^{(k)}_l)}{w^{(k)}_l-z}  \right)\di{\xi} \right|\\
&\leq \sum_{l=1}^M  \frac{\delta}{M}  \max_{\xi \in  \tau_{k,l}}  \left|\frac{f(\xi)}{\xi-z}  -    \frac{f(w^{(k)}_l)}{w^{(k)}_l-z} \right| .
\end{align*}
As the function $\xi \mapsto  \tfrac{f(\xi)}{\xi-z}$ is Lipschitz continuous on $ \tau${,} one has
\begin{align*}
\left| \int_{\tau_k} \frac{f(\xi)}{\xi-z}\di{\xi}  -  \sum_{l=1}^{M}  \frac{c^{(k)}_l}{z-w^{(k)}_l} \right| &\leq 
\sum_{l=1}^M  \frac{\delta}{M} L\, \max_{\xi \in \tau_{k,l}}  \left|\xi - w^{(k)}_l \right|
\leq \frac{\delta^2L}{M} 
\end{align*}
for every $z \in K$. Thus, for all $z\in K${,} one has
\begin{align*}
|f(z) - r(z)| \leq \frac{1}{2\pi} \sum_{k=1}^N \left|  \int_{\tau_k} \frac{f(\xi)}{\xi-z} \di{\xi} -  \sum_{l=1}^{M}  \frac{c^{(k)}_l}{w^{(k)}_l-z} \right| 
\leq  \frac{N}{2\pi}\cdot \frac{\delta^2L}{M} \,
< \frac{\e}{3}.
\end{align*}

$\,$ \hfill {\small $\Box$} 

\medskip
\noindent
{\bfseries Step~3:  Pole shifting} \emph{
Let $\eta = \max_{z \in K} |z|$ and let $b \in \C \setminus K$ with $|b|>2\eta$, $\delta_{kl}:= | w^{(k)}_l -b|$ and $\alpha:=\min_{z\in K} |z-b|$. Then, the rational function
\begin{align*}
r_b(z)= \frac{1}{2\pi i} \sum_{k=1}^{N} \sum_{l=1}^{M} c^{(k)}_l q_{kl}\left(\frac{1}{z-b} \right),
\end{align*}
defined by the polynomials $q_{kl}(z) =z+ (w^{(k)}_l-b) z^2+ \cdots + (w^{(k)}_l -b)^{m_{kl}} z^{m_{kl}+1}$ of degree
 \begin{align}\label{eq:degree_r_b}
  m_{kl}+1 \geq  \frac{\log \left( \frac{\e}{3} \frac{2\pi(\alpha-\delta_{kl})}{ NM\,|c^{(k)}_l|} \right) }    {\log(\delta_{kl}) - \log(\alpha) }
 \end{align}
 approximates $r$ uniformly on $K$, i.e. $r_b$ satisfies
\begin{align*}
 \sup_{z\in K} |  r(z)-r_b(z)|< \tfrac{\e}{3}.
\end{align*}
}
\medskip
\noindent
\emph{Proof of Step~3:}
Since
\begin{align}\label{eq:krantz-geometric-series}
 \frac{1}{z-w^{(k)}_l} = \frac{1}{z-b} \cdot \frac{1}{1-\tfrac{w^{(k)}_l-b}{z-b}}= \frac{1}{z-b} \sum_{v=0}^\infty \left(\tfrac{w^{(k)}_l-b}{z-b}\right)^v
\end{align}
for all $z \in K$, one has
\begin{align*}
|r(z)- r_b(z)| &\leq   \frac{1}{2\pi} \sum_{k=1}^{N} \sum_{l=1}^{M}   | c^{(k)}_l| \,\left|   \frac{1}{w^{(k)}_l-z}   - q_{kl}\left(\frac{1}{z-b} \right) \right| \\
&= \frac{1}{2\pi} \sum_{k=1}^{N} \sum_{l=1}^{M}   | c^{(k)}_l| \,\left|   \frac{1}{w^{(k)}_l-z}   - \frac{1}{z-b} \sum_{p=0}^{m_{kl}} \left(\frac{w^{(k)}_l -b}{z-b} \right)^p \right| \\
&\leq \frac{1}{2\pi} \sum_{k=1}^{N} \sum_{l=1}^{M}    \, \frac{| c^{(k)}_l| }{|z-b|} \sum_{p=m_{kl}+1}^\infty \left|\frac{w^{(k)}_l -b}{z-b} \right|^p  \\
&\leq \frac{1}{2\pi} \sum_{k=1}^{N} \sum_{l=1}^{M}    \, \frac{| c^{(k)}_l| }{\alpha} \sum_{p=m_{kl}+1}^\infty \left|\frac{\delta_{kl}}{\alpha} \right|^p  \\
&\leq \frac{1}{2\pi} \sum_{k=1}^{N} \sum_{l=1}^{M}    \, \frac{| c^{(k)}_l| }{\alpha - \delta_{kl}}  \left(\frac{\delta_{kl}}{\alpha}\right)^{m_{kl}+1}.
\end{align*}
The assertion then follows from observing that \eqref{eq:degree_r_b} is equivalent to 
 \begin{align*}
 \left(\frac{\delta_{kl}}{\alpha}\right)^{m_{kl}+1} \leq \frac{\e}{3} \, \frac{2\pi(\alpha-\delta_{kl}) }{ NM\,|c^{(k)}_l|}.
 \end{align*}
$\,$ \hfill {\small $\Box$}

\medskip
\noindent
{\bfseries Step~4:  Polynomial approximation} \emph{
The polynomial
\begin{align}\label{eq:runge_poly}
p(z):= \frac{1}{2\pi i} \sum_{k=1}^{N} \sum_{l=1}^{M} c^{(k)}_l  \sum_{v=0}^{m_{kl}} (w^{(k)}_l-b)^v p_{kl}^{(v)}(z),
\end{align}
where
\begin{align*}
p_{kl}^{(v)}(z)&:= a^{(klv)}_0 + a^{(klv)}_1z+\cdots + a^{(klv)}_{m^{(v)}_{kl}} z^{m^{(v)}_{kl}+1},\\
a^{(klv)}_s &:= \frac{1}{2\pi i} \int_{\partial B_{r}(0)} \frac{\di{\xi}}{\xi^{s+1} (\xi-b)^{v+1}}, \qquad s=0,1,...,m_{kl}^{(v)}, \quad r:=\tfrac{|b|}{2},
\end{align*}
and
\begin{align}\label{eq:degree_p}
m^{(v)}_{kl} +1 > 
\frac{\log \left( 1+\frac{\e}{3} \cdot \frac{2\pi}{ NM\,r\, |c^{(k)}_l|} \cdot (r-\eta) \cdot (\delta_{kl}-r)\right) } {\log \delta_{kl} - \log r } ,\qquad \eta := \max_{z \in K}|z|{,}
\end{align}
approximates $r_b$ uniformly on $K$, i.e.{,} $p$ satisfies
\begin{align*}
 \sup_{z\in K} |  r_b(z)-p(z)|< \tfrac{\e}{3}.
\end{align*}
}
\medskip
\noindent
\emph{Proof of Step~4:} By construction, for $z \in K$, we have 
\begin{align*}
|  r_b(z)-p(z)| &\leq  \frac{1}{2\pi} \sum_{k=1}^{N} \sum_{l=1}^{M} |c^{(k)}_l|  \left| \sum_{v=0}^{m_{kl}}  \frac{(w^{(k)}_l-b)^v}{(z-b)^{v+1}} - (w^{(k)}_l-b)^v p_{kl}^{(v)}(z) \right|\\
&\leq \frac{1}{2\pi} \sum_{k=1}^{N} \sum_{l=1}^{M} |c^{(k)}_l|  \sum_{v=0}^{m_{kl}}  | w^{(k)}_l-b|^v \, \left| \frac{1}{(z-b)^{v+1}} - p_{kl}^{(v)}(z) \right|.
\end{align*}
Since the function $\frac{1}{(z-b)^{v+1}}$ is holomorphic on $B_r(0)${,} its series expansion around zero is given by the coefficients
\begin{align*}
a^{(klv)}_\mu := \frac{1}{2\pi i} \int_{\partial B_r(0)} \frac{\di{\xi}}{\xi^{\mu+1} (\xi-b)^{v+1}}, \qquad \mu=0,1,2,3,...., \quad r:=\tfrac{|b|}{2}>1.
\end{align*}
That is, 
\begin{align*}
\frac{1}{(z-b)^{v+1}} = \sum_{\mu=0}^\infty a_\mu^{(klv)} z^\mu{,}
\end{align*}
and one has
\begin{align*}
\left| \frac{1}{(z-b)^{v+1}} - p_{kl}^{(v)}(z) \right|\leq \sum_{\mu=m_{kl}^{(v)}+1}^\infty |a_\mu^{(klv)}| |z|^\mu.
\end{align*}
By the Cauchy estimates \cite[Theorem~3.4.1]{greene_krantz_1997}{,} one has
\begin{align*}
|a^{(klv)}_\mu | \leq   \frac{1}{r^{\mu}} \max_{|\xi|=r} \frac{1}{|\xi-b|^{v+1}} =   \frac{1}{r^{\mu}} \cdot \frac{1}{\min_{|\xi|=r} |\xi-b|^{v+1}} = \frac{1}{r^{\mu}} \cdot \frac{1}{ r^{v+1}} .
\end{align*}
Consequently, with $| w^{(k)}_l-b|^v = (\delta_{kl})^v$ for $z \in K$, we have 
\begin{align*}
|  r_b(z)-p(z)| 
&\leq \frac{1}{2\pi} \sum_{k=1}^{N} \sum_{l=1}^{M} |c^{(k)}_l| \,  \sum_{v=0}^{m_{kl}}   \, \frac{(\delta_{kl})^v}{r^{v+1}}  \sum_{\mu=m_{kl}^{(v)}+1}^\infty  \left(\tfrac{|z|}{r}\right)^\mu \\
&\leq \frac{1}{2\pi} \sum_{k=1}^{N} \sum_{l=1}^{M} |c^{(k)}_l| \,  \sum_{v=0}^{m_{kl}}   \, \frac{(\delta_{kl})^v}{r^{v+1}}  \sum_{\mu=m_{kl}^{(v)}+1}^\infty  \left(\tfrac{\eta}{r}\right)^\mu \\
&= \frac{1}{2\pi} \sum_{k=1}^{N} \sum_{l=1}^{M}  |c^{(k)}_l|\, \frac{\left(\tfrac{\eta}{r}\right)^{m_{kl}^{(v)}+1}}{r-\eta}  \sum_{v=0}^{m_{kl}}   \, \left(\frac{\delta_{kl}}{r} \right)^{v} \\
&= \frac{1}{2\pi} \sum_{k=1}^{N} \sum_{l=1}^{M}  |c^{(k)}_l|\, \frac{\left(\tfrac{\eta}{r}\right)^{m_{kl}^{(v)}+1}}{r-\eta}    \, 
 \left( \frac{ 1- \left( \tfrac{\delta_{kl}}{r}\right)^{m_{kl}+1} }{1-\frac{\delta_{kl}}{r}  }     \right) \\
 &< \frac{1}{2\pi} \sum_{k=1}^{N} \sum_{l=1}^{M}  |c^{(k)}_l|\, \frac{r}{(r-\eta)(r-\delta_{kl})}    \, 
 \left( 1- \left( \tfrac{\delta_{kl}}{r}\right)^{m_{kl}+1}      \right) ,
\end{align*}
where in the last inequality we used that $ \left(\tfrac{\eta}{r}\right)^{m_{kl}^{(v)}+1}<1$ since $\eta < r$. By construction{,} it holds $\delta_{kl}>r$ and thus \eqref{eq:degree_p} is equivalent to
 \begin{align*}
 1 - \left(\frac{\delta_{kl}}{r}\right)^{m_{kl}^{(v)}+1} < \frac{\e}{3} \cdot \frac{2\pi}{NM\,r\,|c^{(k)}_l|} \cdot (r-\eta) \cdot (r-\delta_{kl}) .
 \end{align*}
This completes step~4.
$\,$ \hfill {\small $\Box$}

\medskip

The above construction steps yield for every $z\in K$
\begin{align*}
  |f(z) - p(z)| \leq  |f(z) - r(z)| + |r(z) - r_b(z)| + |r_b(z) - p(z)| < \e,
\end{align*}
which completes the proof. \hfill {$\blacksquare$}
\end{proof}

\begin{remark}
\begin{enumerate}
\item[(a)] The set $\bigcup_{k=1}^N \tau_k$ is determined only by the sets $\Omega$ and $K$, and the grid constructed in Step~1. The set $\bigcup_{k=1}^N \tau_k$ also contains all the poles of $r$ and is independent from the quality of approximation. 

\item[(b)] {As $\e>0$ is decreasing the number of poles of the approximating rational function is increasing. The location of the poles is not affected, i.e.,} the poles will still be located on  $\bigcup_{k=1}^N \tau_k$.

\item[(c)] It is shown in \cite[Chapter~12, \S~4]{remmert2013classical} that the line segments $\tau_1,...,\tau_N$ {define a} cycle $\tau$ consisting of finitely many grid polygons. In general{,} it is not possible to find a single closed path $\Gamma$ such that
\begin{align*}
 f(z)=   \frac{1}{2\pi i} \int_{\Gamma} \frac{f(\xi)}{\xi-z} \di{\xi}\qquad \forall \, \, z \in K.
\end{align*}
A counterexample is given \cite[Chapter~12, \S~1]{remmert2013classical}. However, under suitable assumptions on the compact set $K$, the line segments $\tau_1,...,\tau_N$ yield a  single closed grid polygon, cf. \cite[Theorem~9~(b)]{JDE_ensembles_2021}.

\end{enumerate}

\end{remark}

\subsection{Walsh's Theorems}\label{sec:appox_walsh}

In this subsection we shall prove Walsh's refinement of Runge's theorem. We note that Walsh's result is still considerably short of Mergelyan's definite result on complex approximation, cf. \cite[Chapter~III~\S 2~A]{gaier1987}. We note that the proof of Mergelyan's result is not constructive, which is true for almost all results in polynomial approximation.  Recall that, a {\it Jordan arc} is a homeomorphic image of a compact interval. The result is {as} follows, cf. \cite[Bemerkung~$6^{\circ}$]{walsh1926_2}.

\begin{theorem}[Walsh (1926)]\label{thm:walsh-arc}
Let $\gamma$ be a Jordan arc. If $f$ is continuous on $\gamma${,} it can {be} uniformly approximated on $ \gamma$ by a polynomial in $z$.
\end{theorem}

 A proof of Theorem~\ref{thm:walsh-arc} is given at the end of this section.  We emphasize that the following exposition is more explicit in terms of constructing the polynomial {than} the ones in the literature. The significance of the subsequent exposition is that the construction of the polynomial $p$ is traced back to Runge's result by using suitable conformal mappings. Numerical procedures for conformal mappings {are} an active research area and beyond the scope of this paper. At the end of this section we provide more detailed comments and references.

Before we come to the proof{,} we need some preparation. To do so, we have to consider sequences of Jordan curves, or more precisely we provide a notion of convergence for sequences of Jordan curves. Recall that, a {\it Jordan curve} {is a homeomorphic image (within $\C$) of the unit circle $\partial \mathbb{D}$}.

Following Courant \cite{courant1922bemerkungen}, a sequence $(\Gamma_n)_{n \in \N} \subset \C$ of Jordan curves is said to converge to a Jordan curve $\Gamma$ if (i) any accumulation point of any sequence $({z}_n)_{n \in \N}$ with ${z}_n \in \inte \Gamma_n $ lies on $\Gamma$ and (ii) for every $\e>0$ there are $n(\e) \in \N$ and $\delta(\e)>0$ with $\lim_{\e \to 0} \delta(\e) = 0$ such that for every ${z} \in \Gamma$ and for every ${z}_1,{z}_2 \in \inte \Gamma_n \cap B_{\e}(p)$, $n > n(\e)$ there is a polygon $\tau_n$ in $\Gamma_n$ connecting ${z}_1$ and ${z}_2$ so that $\operatorname{diam} \tau_n:= \sup \{|{v}-w|\, {:\, v}, w \in \tau_n\}< \delta(\e)$.

Let $\Omega$ be a proper subdomain of $\C$. Recall that a function $f\colon \Omega \to \D$ is called conformal{,} if it is holomorphic, one-to-one and onto. The following result is due to Courant, cf. \cite{courant1914ueber,courant1922bemerkungen}.

\begin{theorem}[Courant (1914)]\label{thm:courant}
Let $\Omega$ and $(\Omega_n)_{n\in \N}$ be domains containing the origin such that the boundaries $\Gamma=\partial \Omega$ and $\Gamma_n=\partial \Omega_n$ are Jordan curves. Let $\psi \colon \overline{\D} \to \overline{\Omega}$ and $\psi_n\colon \overline{\D} \to \overline{\Omega_n}$, $n\in \N${,} denote continuous functions that are one-to-one, onto, and conformal on $\D$ such that $\psi(0)=\psi_n(0)=0$ and $\psi'(0)>0 $ and $\psi'_n(0)>0$. {Then} $(\psi_n)_{n \in \N}$ converges uniformly on $\overline{\D}$ to $\psi$ if and only if $(\Gamma_n)_{n\in \N}$ converges to $\Gamma$.
 \end{theorem}

Using a different method, Courant's Theorem was also derived by Rad\'{o}  \cite{rado1923representation}. We note that Rad\'{o} used  the Fr\'{e}chet-distance to characterize convergence for sequences of Jordan curves. That these notions are in fact equivalent was shown by Markushevich. A nice and in-depth exposition for this, as well as extensions to more general domains can be found in \cite{Gaier_math_z_1956}.

For a given Jordan curve $\Gamma${,} a sequence of Jordan curves $(\Gamma_n)_{n \in \N}$ that converges to $\Gamma$ from outside, i.e.{,} $\Gamma_n \subset \C \setminus \overline{\inte \Gamma}$, can be constructed via the grid construction in Step~1 in the proof Runge's Theorem. Indeed, as illustrated in Step~1 in the proof Runge's Theorem, it is sufficient to choose a grid consisting of horizontal and vertical lines defined by
\begin{align*}
 x = \frac{k}{2^n} \qquad  y= \frac{l}{2^n}, \quad k,l \in \mathbb{Z}
\end{align*}
and define $\Gamma_n$ as the boundary of the boxes given by the vertices
\begin{align*}
\left(\frac{k}{2^n}, \frac{l}{2^n} \right), \, \, \left(\frac{k+1}{2^n}, \frac{l}{2^n} \right) , \, \, \left(\frac{k+1}{2^n}, \frac{l+1}{2^n} \right), \, \,\left(\frac{k}{2^n}, \frac{l+1}{2^n} \right)
\end{align*}
that intersect with $\overline{\inte \Gamma}$. Next we recap a proof for a refinement of Runge's {theorem} which is also due to Walsh, cf. \cite[Satz]{walsh1926_1}.

\begin{theorem}[Walsh (1926)]\label{thm:walsh_jordan_curve_poly}
Let $\Gamma$ be a Jordan curve and suppose that $f$ is holomorphic on $\inte \Gamma$
and  continuous on $\overline{ \inte \Gamma}$. Then, $f$ can be uniformly approximated on $\overline{ \inte \Gamma}$ by polynomials. That is, for every $\e>0$ there is a polynomial $p$ such that
 \begin{align*}
\sup_{z\in \overline{\inte \Gamma} } |  f (z)-  p(z)| < \e.
 \end{align*}
\end{theorem}
\begin{proof}
 Let $\e>0$ and let w.l.o.g. $z=0 \in \inte  \Gamma$. Using the grid construction as in the proof of Runge's theorem there is {a} sequence $(\Gamma_n)_{n \in \N}$ of grid polygons such that $ \Gamma_n \subset \C \setminus \overline{\inte \Gamma}$ and $(\Gamma_n)_{n \in \N}$ converges to $\Gamma$. {Then} we consider the conformal mappings
\begin{align*}
\Phi &\colon {\inte \Gamma}  \to \D, \quad \,\,\,\, \, \Phi(0)=0, \, \Phi'(0)>0 \\
\Phi_n &\colon {\inte \Gamma_n}  \to \D, \quad \Phi_n(0)=0, \, \Phi_n'(0)>0 .
\end{align*}

\medskip
{\it Step~1: ({Approximation of} $f$ by a sequence of holomorphic functions on $ \overline{\inte \Gamma}$)}. \\
We consider the holomorphic functions
\begin{align*}
g_n &{:=} f \circ \Phi^{-1}\circ \Phi_n \colon {\inte \Gamma_n}  \to \C
\end{align*}
and shall show that for every $\e>0$ there is an $N\in \N$ such that 
\begin{align}\label{eq:g_N-f}
|g_N(z) -  f(z)| < \tfrac{\e}{2} \qquad \text{ for all } \quad z \in \overline{\inte \Gamma}.
\end{align}

Indeed, for $n \in \N$ we consider the Jordan curve $\gamma_n:=\Phi_n(\Gamma) \subset \D$ and note that for every $z \in \overline{\inte \Gamma}$ there exists an $w\in \overline{\inte \gamma_n}$ so that $z=\Phi_n^{-1}(w)$. Thus, together with the maximum modulus theorem~\cite[Corollary~5.3.4]{greene_krantz_1997}, it is sufficient to show that there is a $K=K(\e) \in \N$ such that
\begin{align*}
\sup_{w \in  \gamma_K } |f( \Phi^{-1}(w))- f( \Phi_K^{-1}(w))| < \tfrac{\e}{2}.
\end{align*}

Note that, $\overline{\inte \Gamma}$ is compact and $f$ is continuous, and so $f$ is uniformly continuous on  $\overline{\inte \Gamma}$. 
{Thus} there is a $\delta(\e)>0$ such that $|f(z)-f(w)|< \tfrac{\e}{2}$ for all $z,w \in \overline{\inte \Gamma}$ with $|z-w|<\delta(\e)$. Since the Jordan curves $(\Gamma_n)_{n\in \N}$ converge to $\Gamma$ we can apply Courant's Theorem~\ref{thm:courant} to the sequence $(\Phi^{-1}_n)_{n\in \N}$ and conclude that there is an $N=N(\e) \in \N $ such that
\begin{align*}
|\Phi^{-1}(w) - \Phi^{-1}_N(w)| < \delta(\e)   \qquad \text{ for all } \quad  w \in \gamma_N \subset \D,
\end{align*}
which shows \eqref{eq:g_N-f}.

\medskip
{\it Step~2: ({Approximation of} $g_N$ via Runge's Theorem)}. The function
\begin{align*}
g_N &= f \circ \Phi^{-1}\circ \Phi_N \colon \overline{\inte \Gamma}  \to \C
\end{align*}
is holomorphic on the compact set $\overline{\inte \Gamma}$. Then, by Runge's Theorem~\ref{thm:runge}, for $\tfrac{\e}{2}>0$ there is a polynomial $p$ such that
\begin{align*}
\sup_{z \in \overline{\inte \Gamma} }  |g_N(z)- p(z)| < \tfrac{\e}{2}.
\end{align*}
Consequently, for every $z \in \overline{\inte \Gamma}  $ it follows that
\begin{align*}
|f(z)-p(z)|\leq |f(z)-g_N(z)| +  |g_N(z)- p(z)| < \e.
\end{align*}
 $\hfill \blacksquare$
\end{proof}

\medskip
We will use the latter statement to prove Walsh's extension of the second Weierstrass Theorem to arbitrary Jordan curves, cf. \cite[Satz~I]{walsh1926_2}.

\begin{theorem}[Walsh (1926)]\label{thm:satz_I}
Let $\Gamma$ be a Jordan curve containing the origin in its interior and let the function $f$ be continuous on $\Gamma$. {Then} $f$ can be uniformly approximated on $\Gamma$ by polynomials in $z$ and  $\tfrac{1}{z}$.
\end{theorem}

\begin{proof}
Let $\e>0$. By Caratheodory's Theorem~\cite[Theorem~13.2.3]{greene_krantz_1997}{,} there is a continuous and bijective map $\Phi\colon \overline{\inte \Gamma} \to \overline{\D} $ so that $\Phi(0)=0$ which is conformal in $\D$. Then, the mapping $$g:=f\circ \Phi^{-1} \colon \partial \D \to \C$$ is continuous. Consider the sequence of Fej\'{e}r polynomials $(F_{g,n})_{n\in \N}\colon \partial \D\to \C$ defined by
\begin{align*}
F_{g,n}(w) :=   \sum_{k=1}^{n-1} \frac{n-|k|}{n} \hat g (-k)\, \,  (\tfrac{1}{w})^{k}  + \sum_{k=0}^{n-1} \frac{n-|k|}{n} \hat g (k)\, \,  w^{k},
\end{align*}
where the $\hat g(k)$ denote the Fourier coefficients of $g$, cf.~\eqref{def:poly_Fejer_coeff}.
Since for every $z \in \Gamma$ there is a unique $w \in \partial \D$ such that $z= \Phi^{-1}(w)$, by Theorem~\ref{thm:weier-circ}~(a), there is a $N=N(\e) \in \N$ such that
\begin{align*}
 \sup_{z \in \Gamma}|f(z)- F_{g,N}(\Phi(z))| =   \sup_{w\in \partial \D} |g(w) -F_{g,N}(w)| < \frac{\e}{3}.
\end{align*}
To show the claim, we shall prove that {on $\Gamma$, the function $\Phi$ can be uniformly approximated on $\Gamma$ by polynomials in $z$ and the function $\tfrac{1}{\Phi}$ can be uniformly approximated  by polynomials in $z$ and $\tfrac{1}{z}$.}

Let $\delta>0$; we will specify it later, cf.~\eqref{eq:walsh-circ-deltas}. By Theorem~\ref{thm:walsh_jordan_curve_poly}, there is a polynomial $p_1$ such that
\begin{align*}
| {\Phi(z)} - p_1(z)| < \delta \qquad \forall \, \, z \in \overline{\inte \Gamma}.
\end{align*}
Next, we show that there is a polynomial $p_2$ such that
\begin{align*}
\left| \frac{1}{\Phi(z)}- \frac{a_{-1}}{z} + p_2(z)\right| < \delta  \qquad \text{ for all } \, \, z \in \Gamma,
\end{align*}
where $a_{-1}$ is the residue of $\tfrac{1}{\Phi}$ at zero, i.e.
$$a_{-1} = \frac{1}{2\pi i} \int_{\partial B_r(0)} \frac{1}{\Phi(\xi)} \di{\xi}, \qquad r >0 \, \, \text{such that } \overline{B_r(0)} \subset \inte \Gamma. $$

To see this, recall that $\Phi$ is conformal on $\inte \Gamma$ with $\Phi(0)=0$. Hence, the function $\tfrac{1}{\Phi}$ is meromorphic on $\inte \Gamma$ with a simple pole at $0$. Moreover, the function $ z \mapsto \tfrac{1}{\Phi(z)} - \tfrac{a_{-1}}{z} $ has a removable singularity at zero. Let $\Psi\colon \inte \Gamma\to \C$ denote its holomorphic extension, i.e.
\begin{align}\label{eq:def-phi_e-walsh-curve}
\Psi(z) = \begin{cases} 
\frac{1}{\Phi(z)} -  \frac{a_{-1}}{z}  & z \in \inte \Gamma\setminus \{0\}\\
a_0= \frac{1}{2\pi i} \int_{\partial B_r(0)} \frac{1}{\xi \Phi(\xi)} \di{\xi} & z=0,
\end{cases} 
\end{align}
which is continuous on $\overline{ \inte \Gamma}$. Then, by Theorem~\ref{thm:walsh_jordan_curve_poly}{,} there is {a} polynomial $p_2$ such that
\begin{align*}
\left| \Psi(z)- p_2(z)\right| < \delta  \qquad \text{ for all } \, \, z \in \overline{ \inte \Gamma }.
\end{align*}
Next, we shall show that{,} if we choose $\delta>0$ such that
\begin{align}\label{eq:walsh-circ-deltas}
   \sum_{k=0}^{N-1} \left( (1+\delta)^k -1 \right) < \frac{\e}{3\,c}, \qquad c := \max_{k=-N+1,...,N-1}  {\Big|}\tfrac{N-|k|}{N} \hat g(k){\Big|}
\end{align}
one has
\begin{align}\label{eq:walsh-jordan-kniff}
 \left| \sum_{k=0}^{N-1} \frac{N-|k|}{N} \hat g (k)\, \,  \Phi(z)^k - \sum_{k=0}^{N-1} \frac{N-|k|}{N} \hat g (k)\, \,  p_1(z)^k\right| < \frac{\e}{3} \qquad \forall \, \, z\in \Gamma.
\end{align}

To see this, note that $|\Phi(z)| =1$ for all $z \in \Gamma$ and, thus, 
\begin{align*}
  |p_1(z)| < 1+ \delta  \qquad \text{ for all } z \in \Gamma.
\end{align*}
Moreover, we recall that for all $v,w \in \C$ it holds
\begin{align*}
  |v^n-w^n| \leq |v-w| \sum_{k=0}^{n-1} |v|^k \, |w|^{n-1-k}.
\end{align*}
The combination of the latter yields that for all $z\in \Gamma$ one has
\begin{multline*}
  \left| \sum_{k=0}^{N-1} \frac{N-|k|}{N} \hat g (k)\, \,  \Phi(z)^k - \sum_{k=0}^{N-1} \frac{N-|k|}{N} \hat g (k)\, \,  p_1(z)^k \right| \\ 
 \leq \sum_{k=0}^{N-1}  \left|\frac{N-|k|}{N} \hat g (k)\right|  \left| \Phi(z)^k - p_1(z)^k\right| \\
 \leq \sum_{k=0}^{N-1}  \left|\frac{N-|k|}{N} \hat g (k)\right|  \left| \Phi(z) - p_1(z)\right| \,\,
 \sum_{l=0}^{k-1}  \,| p_1(z)|^l \leq
  c \sum_{k=0}^{N-1} \left\{ (1+\delta)^k -1 \right\}< \frac{\e}{3}.
\end{multline*}

Similarly, for all $z\in \Gamma$ we get
\begin{align*}
 \left| 
 \sum_{k=1}^{N-1} \frac{N-|k|}{N} \hat g (-k)\, \,  (\tfrac{1}{\Phi(z)})^k - \sum_{k=1}^{N-1} \frac{N-|k|}{N} \hat g (-k)\, \,  \left(\tfrac{a_{-1}}{z} + p_2(z)\right)^k
 \right| 
 < \frac{\e}{3}.
\end{align*}
Finally, defining the polynomial 
\begin{align*}
P(z,\tfrac{1}{z}) 
= \sum_{k=1}^{N-1} \frac{N-|k|}{N} \hat g (-k)\, \, \left(\tfrac{a_{-1}}{z} + p_2(z)\right) ^k + \sum_{k=0}^{N-1} \frac{N-|k|}{N} \hat g (k)\, p_1(z)^k{,}
\end{align*}
we derive the following estimate
\begin{align*}
\left|f(z) - P(z,\tfrac{1}{z}) \right| &\leq \left|f(z) - F_{g,N}(\Phi(z)) \right| \\
&\quad+ \left|    \sum_{k=1}^{N-1} \frac{N-|k|}{N} \hat g (-k)\,\left( \, \left( \tfrac{1}{\Phi(z)} \right)^k - 
\, \, \left(\tfrac{a_{-1}}{z} + p_2(z)\right) ^k \right)
 \right| \\
 &\quad+ \left|\sum_{k=0}^{N-1} \frac{N-|k|}{N} \hat g (k)\left(\Phi(z)^k - p_1(z)^k \right) \right| < 3\, \tfrac{\e}{3} = \e.
\end{align*}
This shows the assertion. $\hfill \blacksquare$
\end{proof}

\medskip

\noindent
{\it Proof of Theorem~\ref{thm:walsh-arc} } 
Without loss of generality we assume that $\gamma$ does not contain the origin\footnote{This can be seen by considering a conformal mapping $z \mapsto vz+w$ for some $v \in  \C \setminus \{0\}$ and $w \in \C$.}. Then, the Jordan arc can be extended to a Jordan curve $\Gamma$ so that $0 \in \inte \Gamma$. Also the function $f$ can be continuously extended to $\Gamma$. This can be {achieved} by defining
\begin{align}\label{eq:extension_f_arc_to_curve}
f(z)= w_1 + (w_2-w_1)\frac{z-z_1}{z_2-z_1} \qquad \forall \, \, z \in \Gamma \setminus \gamma,
\end{align}
where $w_1$ and $w_2$ are the values of $f$ at the end-points $z_1$ and $z_2$ of $\gamma$, respectively. By Caratheodory's Theorem~\cite[Theorem~13.2.3]{greene_krantz_1997} there is a continuous and bijective map $\Phi\colon \overline{\inte \Gamma} \to \overline{\D} $ so that $\Phi(0)=0$ which is conformal on $\inte \Gamma$.

Let $\e>0$. By Theorem~\ref{thm:satz_I}, for the continuous function $f \colon \Gamma \to \C$ there are $N=N(\e) \in \N$ and $a_{-1} \in \C$ {as well as} polynomials $p_1,p_2${,} such that the polynomial
\begin{align*}
P_{f\circ \Phi^{-1},N}(z,\tfrac{1}{z}) &:= 
\sum_{k=1}^{N-1}  c_{-k}\, \, \left(\tfrac{a_{-1}}{z} + p_2(z)\right) ^k + \sum_{k=0}^{N-1}  c_k\, p_1(z)^k
\end{align*}
in $z$ and $\tfrac{1}{z}$ satisfies
\begin{align*}
\sup_{z \in {\Gamma}} |f(z) -P_{f\circ\Phi^{-1},N}(z,\tfrac{1}{z}) |< \frac{\e}{2}.
\end{align*}
Here{,} the coefficients $c_{-N+1},...,c_{N-1}$ are given by the Fourier coefficients corresponding to the function $f \circ \Phi^{-1}:\partial \D \to \C$, cf. \eqref{def:poly_Fejer_coeff}, i.e.
\begin{align*}
c_k:=  \tfrac{N-|k|}{N} \widehat{f\circ \Phi^{-1}} (k) , \qquad k \in \{-N+1,...,N-1\}. 
\end{align*}

Furthermore, on $\gamma$ the function $z \mapsto \tfrac{1}{z}$ can be uniformly approximated by polynomials in $z$. To see this, let $\gamma_{0,\infty}$ be a Jordan arc connecting the origin and $\infty$ such that $ \gamma_{0,\infty} \cap  \gamma = \emptyset$. Then,  $\Omega:= \C \setminus  \gamma_{0,\infty}$ is simply connected, $\infty \not \in \inte \Omega$ and the function $z \mapsto \tfrac{1}{z}$ is holomorphic on $\inte \Omega$.  Let 
\begin{align*}
 \nu:= \max_{z \in \gamma} \left|\tfrac{1}{z} \right|, \, \,  
 \mu:= \max_{z \in \gamma} \left|p_2(z) \right|
  \quad \text{and} \quad c:= \max_{k=-N+1,....,-1} \left|  c_k \right|. 
\end{align*}
Moreover, choose $\delta>0$ such that
\begin{align}\label{eq:delta-walsh_arc}
\delta  \cdot \frac{1-(2\delta + |a_{-1}|\nu +\mu)^{N-1}}{1-(2\delta + |a_{-1}|\nu +\mu)}   < \frac{\e}{2 c(N-1)}.
\end{align}
By Runge's Theorem~\ref{thm:runge}{,} there is {a} polynomial $p_3$ in $z$ such that
\begin{align*}
\sup_{z \in \gamma} |\tfrac{a_{-1}}{z} -p_3(z)|< \delta.
\end{align*}
Next, we define the polynomial
\begin{align}\label{eq:def-poly-walsh-arc}
p(z) := \sum_{k=0}^{N-1}  c_k\,p_1(z)^k +    \sum_{k=1}^{N-1}  c_{-k} \, \big( p_3(z) + p_2(z)\big)^k.
\end{align}
Similar to the proof of Theorem~\ref{thm:satz_I}{,} we get
\begin{align*}
\left|  \sum_{k=1}^{N-1}  c_{-k}\,\big( \tfrac{a_{-1}}{z}+ p_2(z)\big)^k  -  \sum_{k=1}^{N-1}  c_{-k}\, \big(p_3(z) + p_2(z)\big)^k \right| 
\\ \leq 
c \sum_{k=1}^{N-1} \left| ( \tfrac{a_{-1}}{z}+ p_2(z))^k  -   (p_3(z) + p_2(z))^k  \right|
\\ \leq 
c\, \delta  \sum_{k=1}^{N-1}   \sum_{l=0}^{k-1} \left|  \tfrac{a_{-1}}{z}+ p_2(z)\right|^l  \left| p_3(z) + p_2(z)  \right|^{k-1-l}.
\end{align*}
Using the {triangle} inequality, for every $z \in \gamma${,} one has
\begin{align*}
 | p_3(z)|< \delta + |\tfrac{a_{-1}}{z}| \leq \delta + |a_{-1}|\nu
\end{align*}
and $ | p_3(z) + p_2(z) |<  \delta + |a_{-1}|\nu + \mu$
as well as
\begin{align*}
 |\tfrac{a_{-1}}{z} + p_2(z) | \leq  |\tfrac{a_{-1}}{z} - p_3(z)| + |  p_3(z) -   p_2(z) |            < 2 \delta + |a_{-1}|\nu + \mu.
\end{align*}
Thus, by \eqref{eq:delta-walsh_arc} we have
\begin{align*}
\left|  \sum_{k=1}^{N-1}  c_{-k}\,\big( \tfrac{a_{-1}}{z}+ p_2(z)\big)^k  -  \sum_{k=1}^{N-1}  c_{-k}\, \big(p_3(z) + p_2(z)\big)^k \right| 
\\ \leq 
c\, \delta \, (N-1) \sum_{k=0}^{N-2}   \left(  2 \delta + |a_{-1}|\nu + \mu  \right)^{k} < \frac{\e}{2}.
\end{align*}
Consequently, it holds
\begin{multline*}
\sup_{z \in  \gamma} |f(z) -p(z)| \leq \sup_{z \in \Gamma} |f(z) -p(z)| 
 \leq 
\sup_{z \in  \Gamma} |f(z) -P_{f\circ \Phi^{-1},N}(z,\tfrac{1}{z})| \\+
\sup_{z \in  \Gamma}
\left|  \sum_{k=1}^{N-1}  c_{-k}\,\big( \tfrac{a_{-1}}{z}+ p_2(z)\big)^k  -  \sum_{k=1}^{N-1}  c_{-k}\, \big(p_3(z) + p_2(z)\big)^k \right| 
<\e.
\end{multline*}
This shows the assertion.
$\hfill \blacksquare$
\medskip

We close this section {with} some comments on the construction of conformal mappings and their numerical computation. There are different approaches {available} in the literature. Peter Henrici presents a method to determine conformal mappings by solving a Dirichlet boundary problem. The solution of it together with its harmonic conjugate determines the conformal map, cf. \cite[Theorem~16.5a]{henrici1993applied}. Moreover, Henrici describes numerical techniques to solve the Dirichlet problem by means of numerical methods to solve Symm's integral equation of first kind, cf.  \cite[\S~16.6~V]{henrici1993applied}.  Besides {that}, Tobin. A. Driscoll and Llyod N. Trefethen depict procedures, as suggested by Peter Henrici, {to compute conformal maps by using the Schwarz-Christoffel mapping}, cf. \cite{trefethen_schwarz}.  Another, recent approach to compute conformal maps is the Zipper-algorithm \cite{marshall2007convergence}. For earlier approaches to this area we refer to \cite{gaier1964konstruktive,gaier1976integralgleichungen}. Furthermore, approximation techniques for conformal mappings based on Bergman kernels are outlined in \cite{gaier1987} and \cite{greene_krantz_1997}. Finally, we mention the monograph \cite{Papamichael_book} {which} contains a profound literature review and many more references.

\section{Discrete-time single-input systems}\label{sec:poly_si_methods}

In this section we consider single-input pairs $(A,b)\in C_{n,n}(\p) \times C_n(\p)${,} where the parameter space $\p$ is a Jordan arc. Especially, for given $f\in C_n(\p)$ and $\e>0${,} we tackle the central problem {of} how to determine suitable inputs such that the origin is steered into the $\e$-neighborhood of $f$. The basic tools for the constructions will be the results of Bernstein, Runge, Weierstrass and Walsh from  the previous section. Indeed, as pointed out in the introduction, the solution to
\begin{align*} \label{eq:def-ensemble}\begin{split}
x _{ t + 1 }( \theta ) = A ( \theta ) x _ t ( \theta ) + b ( \theta ) u _t \end{split}
\end{align*}
can be represented by a polynomial $p(z)= u_0 z^t + u_1z^{t-1} + \cdots +u_{t-1}$, i.e.
\[
 \varphi(t,u,0) = p(A)b.
\]
Hence, given a target function $f \in C_n(\p)$ and $\e>0${,} we shall explore {the} methods of Section~\ref{sec:approx} to derive a suitable polynomial $p$ so that
\begin{align*}
  \| \varphi(t,u,0) -f\|_\infty = \|p(A)b -f\|_\infty < \varepsilon.
\end{align*}

\medskip
\noindent
{\it Starting point: Known sufficient conditions}

As a starting point{,} we take the necessary and sufficient conditions for uniform ensemble reachability developed in \cite[Theorem~4, Corollary~3]{JDE_ensembles_2021}. Since these are essential, we provide a short recap. That is, if the pair $(A,b)$ is uniformly ensemble reachable, the following necessary conditions are satisfied:
\begin{description}
\item[(N1)] The pair $(A(\theta),b(\theta))$ is reachable for every $\theta\in \p$.
\item[(N2)] For any pair of distinct parameters $\theta,\theta'\in \p$, the spectra of $A(\theta)$ and $A(\theta')$ are disjoint:
\begin{equation*}
\sigma (A(\theta))\cap \sigma (A(\theta'))=\emptyset.
\end{equation*}
\end{description}
Further, $(A,b)$ is uniformly ensemble reachable if it satisfies (N1){,} (N2) and one of the following sufficiency conditions:
\begin{description}
\item[(S1)] {For each $\theta \in \p$ the characteristic polynomial of $A(\theta)$ takes} the form
\begin{equation*}
z^n - (a_{n-1}z^{n-1} + \cdots +a_1z + a_0(\theta))
\end{equation*}
for some $a_{n-1},...,a_1 \in \C$ and $a_0 \in C(\p)$.
\item[(S2)] $A(\theta)$ has simple eigenvalues for each $\theta\in \textbf{P}$.
\end{description}

A few comments are in order. First, we note {that} the conditions are valid for continuous-time systems and discrete-time systems. Also, we note that a pointwise characterization for uniform ensemble reachability is currently not available. In \cite[Theorem~3~(d)]{JDE_ensembles_2021} it is {shown that for single input pairs $(A,b)$ a necessary condition for uniform ensemble reachability is that the eigenvalues are simple on an open dense subset of the parameter space $\p$. In this sense condition (S2) is not far from being necessary.}
Furthermore, in certain cases, like the controlled harmonic oscillator, cf. \cite{helmke2014uniform}, it happens that the sufficiency conditions~(S1) and~(S2) are satisfied at the same time. Besides, to the best of the author's knowledge, other general sufficient conditions are not known. Moreover, we note that{,} due to the condition (N2) the function $a_0$ in (S2) is necessarily injective. Thus, $a_0\colon \p \to a_0(\p)$ is one-to-one and onto and so $a_0(\p)$ defines a Jordan arc. It will turn out that in the special case{,} where $a_0(\p)$ is a subset of the real axis, based on Theorem~\ref{thm:gzyl}, a constructive method can be obtained {by} using Bernstein polynomials. The general case is covered by Walsh's Theorem~\ref{thm:walsh-arc}. This method {turns} out to be more involved and will require also the Theorems~\ref{thm:weier-circ},~\ref{thm:runge},~\ref{thm:walsh_jordan_curve_poly} and~\ref{thm:satz_I}.

{Before we begin with the constructions,}
we note that we start the constructions based on the conditions (N1), (N2) in (S1) or ({S}2) and discuss additional required conditions along the lines. The corresponding results are formulated subsequently to each construction. This has the advantage that the results and the corresponding polynomials can be spotted simultaneously.

\subsection{Methods for condition~(S1)}\label{sec:meth_A} 

In this {subsection,} we will explore {the} circumstances such that the verification of
the sufficient condition~(S1) (cf. proof of Theorem~4 in \cite{JDE_ensembles_2021}) gives rise to methods for the construction of suitable $T>0$ and inputs $u=(u_0,...,u_{T-1})$.  {Let's} we assume that the pair $(A,b)$ satisfies the conditions (N1), (N2) and (S1). Also, let $f \in \operatorname{C}_n(\p)$  and $\varepsilon>0$ be given. By condition~(N1), the matrix
\[
 T(\theta) = \begin{pmatrix} b(\theta) & A(\theta)b(\theta) & \cdots & A(\theta)^{n-1}b(\theta)
             \end{pmatrix}
\]
is continuously invertible for every $\theta \in \p$ and one has
\begin{equation*}\label{eq:canonical}
\tilde A(\theta)=
T(\theta)^{-1}A(\theta)T(\theta) = 
\left(\begin{array}{ccccc}
   0 & \hdots & \hdots & 0 & a_0(\theta)\\
   1 & \ddots & & \vdots & a_1\\
   0 &\ddots & \ddots & \vdots & \vdots\\
   \vdots & \ddots & \ddots & 0 & \vdots\\
   0 & \hdots & 0 & 1 & a_{n-1}
\end{array}\right)
,\qquad 
\tilde b= T(\theta)^{-1}b(\theta)=
 \begin{pmatrix}
1\\
0\\
 \vdots \\
0
\end{pmatrix},
 \end{equation*}
where $a_0(\theta),a_1,...,a_{n-1}$ denote the coefficients of the characteristic polynomial of $A(\theta)$. Writing the solution formula in terms of a polynomial $p${,} we get
\begin{align*}
\| \varphi(T,u,0) - f\|_{\infty}  &= \| p(A)b - f\|_{\infty} = \|T ( p(\tilde A)\tilde b - T^{-1}f)\|_{\infty} \\&\leq \|T\|_{M,\infty} \,\| p(\tilde A)\tilde b - T^{-1}f\|_{\infty}, 
\end{align*}
where $\|\cdot \|_M$ {is} a matrix norm that is submultiplicative to $\|\cdot\|$ and
$$\|T\|_{M,\infty} := \sup_{\theta \in \p} \|T(\theta)\|_M.$$
As shown in \cite{JDE_ensembles_2021}, for 
\begin{align}\label{eq:def-p-basic}
p(z)= \sum_{k=1}^n p_k\big(z^n-a_{n-1}z^{n-1}-\cdots - a_1z\big)\, z^{k-1}{,}
\end{align}
we get
\begin{align*}
p(\tilde A(\theta)) \tilde b
 =\sum_{k=1}^n p_k \big(  A(\theta)^n -a_{n-1} A(\theta)^{n-1} -\cdots {-}a_1 A(\theta) \big)A(\theta)^{k-1}\tilde b
 =\begin{pmatrix}
   p_1(a_0(\theta))\\
   \vdots\\
   p_n(a_0(\theta))
  \end{pmatrix}
\end{align*}
and, hence,
\begin{align*}
\| \varphi(T,u,0) - f\|_{\infty}  \leq \|T\|_{M,\infty} \,\sup_{k=1,...,n} \sup_{\theta \in \p} | p_k(a_0(\theta)) - (T^{-1}f)_k(\theta)| ,
\end{align*}
where $(T^{-1}f)_k(\theta)$ denotes the $k$th component of the vector-valued function $\theta \mapsto T^{-1}(\theta)f(\theta)\in \C^n$. Moreover, by (N2) and (S1), the function $\theta \mapsto a_0(\theta)$ is one-to-one. {Consequently, we get}
\begin{align*}
\| \varphi(T,u,0) - f\|_{\infty}  <  \e{,}
\end{align*}
if we can find polynomials $p_1,....p_n$ satisfying
\begin{align}\label{eq:S1_basic}
\sup_{z \in a_0(\p)} | p_k(z) - (T^{-1}f)_k \circ a_0^{-1}(z)|  < \frac{\e}{\|T\|_{M,\infty}} \quad \text{ for all } \, k=1,...,n.
\end{align}
Depending on the structure of $a_0(\p)${,} we will present two constructive methods.

\subsubsection*{Case~1: $a_0(\p)=[a,b]\subset \R$ defines a compact interval}

We assume that for each $k=1,...,n$ the function $z \mapsto (T^{-1}f)_k\circ a_0^{-1}(z)$ is Lipschitz continuous\footnote{This is the case, e.g.{,} if the pair $(A,b)$ is in $C^1_{n,1}(\p)\times C_n^1(\p)$ {and $f \in \operatorname{Lip}_n(\p)$}.} and let $L_k{>}0$ denote the Lipschitz constant. By Theorem~\ref{thm:gzyl}, the condition \eqref{eq:S1_basic} is fulfilled if we take $p_k$ as the Bernstein polynomial, cf. \eqref{def:poly_Bernstein}, corresponding to the function $z \mapsto (T^{-1}f)_k\circ a_0^{-1}(z)$ with degree $n_k \geq 3$ satisfying
\begin{equation*}
\sqrt{2} \left( 4\,M_{(T^{-1}f)_k\circ a_0^{-1}} + \frac{(b-a)\,L_{k}}{2} \right) {\sqrt\frac{\log n_k}{n_k}} < \frac{\e}{\|T\|_{M,\infty}}.
\end{equation*}
Denoting this Bernstein polynomial by $ B_{(T^{-1}f)_k\circ a_0^{-1},n_k}${,} we conclude that for
\begin{align*}
p(z)= \sum_{k=1}^n  B_{(T^{-1}f)_k\circ a_0^{-1},n_k}\big(z^n-a_{n-1}z^{n-1}-\cdots - a_1z\big)\, z^{k-1}
\end{align*}
it holds
\begin{align*}
\| \varphi(T,u,0) - f\|_{\infty} \leq \|T\|_{M,\infty}\, \|p(\tilde A) \tilde b- T^{-1}f\|_\infty <  \e .
\end{align*}
Considering the monomial representation
\begin{align}\label{def:poly_Bernstein_S1}
p(z)= \sum_{k=1}^n  B_{(T^{-1}f)_k\circ a_0^{-1},n_k}\big(z^n-a_{n-1}z^{n-1}-\cdots - a_1z\big)\, z^{k-1} =\sum_{l=0}^N \alpha_l z^l{,}
\end{align}
we get  the following {result.}
\medskip

\begin{theorem}
Let $(A,b) \in C_{n,n}(\p)\times C_n(\p)$ satisfy (N1), (N2) and (S1). Let $\varepsilon>0$ and suppose that $f \in C_n(\p)$ is such that $(T^{-1}f)_k \circ a_0^{-1} \in \operatorname{Lip}(a_0(\p))$ for all $k=1,...,n$ and assume that  $a_0(\p)$ is a compact real interval.  Then, for
\begin{align*}
T=N+1 \qquad \text{and} \qquad  u=(u_0,...,u_{T-1}) = (\alpha_N,...,\alpha_{0})
\end{align*}
defined in \eqref{def:poly_Bernstein_S1}{,} one has
\begin{align*}
\| \varphi(T,u,0) - f\|_{\infty} < \varepsilon.
\end{align*}
In particular, the result is true if $(A,b) \in C^1_{n,n}(\p)\times C^1_n(\p)$ and $f\in \operatorname{Lip}_n(\p)$.
\end{theorem}

Next we consider the more general problem of constructing {polynomials $p_k$} satisfying \eqref{eq:S1_basic} if $a_0(\p)$ defines a Jordan arc. It turns out that this is more involved.  We shall see how the proof of Theorem~\ref{thm:walsh-arc} gives rise to a constructive method.

\subsubsection*{Case~2: $a_0(\p)$ defines a Jordan arc}

Without loss of generality\footnote{The case $0 \in a_0(\p)$ can be reduced to the case $0\not \in a_0(\p)$ by considering a conformal map $z \mapsto vz +w$ for some $v\in \C\setminus \{0\}$ and $w \in \C$.} we assume that $a_0(\p)$ does not contain the origin.
In a first step we extent $a_0(\p)$ to a Jordan curve $\Gamma_0$ such that the origin lies inside the Jordan curve. Then, take a continuous mapping $\Phi_0 : \overline{\inte \Gamma_0} \to \overline{\mathbb D}$ which is conformal on $\inte \Gamma_0$ so that $\Phi_0(0)=0$. The existence of such a mapping is ensured by Caratheodory's Theorem, cf. \cite[Thm.~13.2.3]{greene_krantz_1997}. 

\medskip
\noindent
{\it Step 1.} \textit{Extension of the target function to the Jordan curve}\\
Here{,} we extend the functions  $ z \mapsto (T^{-1}f)_k\circ a_0^{-1}(z)$ to $\Gamma_0$, where $T\in C_{n,n}(\p)$ denotes the changes of coordinates introduced at the beginning of Section~\ref{sec:meth_A}. This can be done as follows. For simplicity{,} we do not introduce a new notation. For all $z \in \Gamma_0 \setminus a_0(\p)$ we define
 \begin{align*}
(T^{-1} f)_k \circ a_0^{-1}(z):= w_{k,1} + (w_{k,2}-w_{k,1})\frac{z-z_{k,1}}{z_{k,2}-z_{k,1}},
\end{align*}
where $w_{k,1}$ and $w_{k,2}$ are the values of $(T^{-1}f)_k\circ a^{-1}_0$ at the end-points $z_{k,1}$ and $z_{k,2}$ of $a_0(\p)$, respectively. As a result{,} we obtain a continuous function
\[
 h_k := (T^{-1}f)_k\circ a_0^{-1} \circ \Phi_0^{-1}  \colon \partial\mathbb{D} \to \mathbb C.
\]
In the following{,} we assume that $a_0 \in C^2(\p)$ and $\p$ has a parametrization which {is} twice continuously differentiable{. Also} we suppose that the extension $\Gamma_0$ has the same properties. Then, by Kellogg's Theorem \cite[Theorem~4.3]{garnett_marshall_2005}, the functions $h_k$ satisfy a Lipschitz condition provided $(T^{-1}f)_k$ is Lipschitz continuous\footnote{This is the case if $(A,b) \in C_{n,n}^1(\p) \times C_n^1(\p)$ {and $f \in \operatorname{Lip}_n(\p)$}.}. We assume that this is the case and let $L_k>0$ denote the corresponding Lipschitz constant.

\medskip
\noindent
{\it Step 2.} \textit{Ansatz for the polynomial $p$}\\
To derive a suitable polynomial $p${,} we take again (cf.~\eqref{eq:def-p-basic})
\begin{align*}
p(z)= \sum_{k=1}^n p_k\big(z^n-a_{n-1}z^{n-1}-\cdots - a_1z\big)\, z^{k-1}.
\end{align*}
Based on the exploration in Section~\ref{sec:appox_walsh} (in particular~\eqref{eq:def-poly-walsh-arc}), we will construct the polynomials $p_1,...,p_n$  as follows. For each $k =1,...,n$ we consider the following Fej\'{e}r{-}like polynomial
\begin{align*}
 p_k(z) = \sum_{l=0}^{n_k-1} c_l \, \big(p_{k,1}(z)\big)^l + \sum_{l=1}^{n_k-1} c_{-l} \,   \big(p_{k,2}(z) + p_{k,3}(z)\big)^l,
\end{align*}
where  
\begin{align*}
 c_l:= \frac{n_k - |l|}{n_k}\,  \widehat{h_k} (l),\qquad \, l = -n_k+1,....,n_k-1
\end{align*}
denote the Fourier coefficients of $h_k$, cf.~\eqref{def:poly_Fejer_coeff}. We proceed with a simple but instructive step. By adding and {subtracting} appropriate terms and using the {triangle} inequality, one has
\begin{align*}
| (T^{-1}f)_k\circ a_0^{-1}(z) &-  p_k(z)| \\
&\leq 
\left|   (T^{-1}f)_k\circ a_0^{-1}  (z) -  \sum_{l=0}^{n_k-1} c_l \, \big(\Phi_0(z)\big)^l  - 
\sum_{l=1}^{n_k-1} c_{-l} \,   \big(\tfrac{1}{\Phi_0(z)}\big)^l\right|\\
&\quad+ \left|\sum_{l=0}^{n_k-1} c_l \, \big(\Phi_0(z)^l- p_{k,1}(z)^l\big)        \right|
\\ &\quad
+ \left|\sum_{l=1}^{n_k-1} c_{-l} \, \Big( \, \big( \tfrac{a_{-1}}{z} +p_{k,2}(z)\big)^l -  \big(\tfrac{1}{\Phi_0(z)}\big)^l \Big)
\right|
\\ &\quad
+ \left|\sum_{l=1}^{n_k-1} c_{-l} \, \Big( \,    \big(p_{k,2}(z) + p_{k,3}(z) \big)^l
- \big( \tfrac{a_{-1}}{z} +p_{k,2}(z)\big)^l  \Big)
\right|,
\end{align*}
where
$$a_{-1} = \frac{1}{2\pi i} \int_{\partial B_r(0)} \frac{1}{\Phi_0(\xi)} \di{\xi}, \qquad r >0 \, \, \text{such that } \overline{B_r(0)} \subset \inte \Gamma_0$$ 
denotes the residue of $\tfrac{1}{\Phi_0}$ at zero.  

\medskip
\noindent
{\it Step 3.} \textit{Derivation of $n_k \in \N$}\\
Applying Weierstrass Second Theorem~\ref{thm:weier-circ}~(b) to the function
\[
 h_k = (T^{-1}f)_k\circ a_0^{-1} \circ \Phi_0^{-1}  \colon \partial\mathbb{D} \to \mathbb C{,}
\]
we get
\begin{align*}
\sup_{w \in \partial \D}\left|h_k (w) -  \sum_{l=0}^{n_k-1} c_l \, w^l  - 
\sum_{l=1}^{n_k-1} c_{-l} \,   \big(\tfrac{1}{w}\big)^l\right| 
< 2\sqrt{2}\pi L_{k} \frac{\ln n_k }{n_k} {.}
\end{align*}
As $\Phi_0$ is one-to-one and onto and taking $n_k \in \N$ so that 
\begin{align*}
  2\sqrt{2}\pi L_{k} \frac{\ln n_k }{n_k}   \leq  \frac{\e}{4\,\|T\|_{M,\infty}},
\end{align*} 
we have
\begin{align*}
\sup_{z \in \Gamma_0} \left| (T^{-1}f)_k\circ a_0^{-1} (z) -  \sum_{l=0}^{n_k-1} c_l \, \big(\Phi_0(z)\big)^l  - 
\sum_{l=1}^{n_k-1} c_{-l} \,   \big(\tfrac{1}{\Phi_0(z)}\big)^l\right| 
< \frac{\e}{4 \, \|T\|_{M,\infty}}. 
\end{align*}

{Hence} it remains to determine the polynomials $p_{k,1},p_{k,2},p_{k,3}$.
This will be carried out in the following three steps. To this end, we pick $\delta_{k}>0$ so that
\begin{align*}
  \sum_{l=0}^{n_k-1} \left( (1+\delta_k)^l -1 \right) < \frac{\e}{4\,c_{\max} \| T\|_{M,\infty}}, \qquad c_{\max} := \max \{{\big|}c_{-n_k+1}{\big|},...,{\big|}c_{n_k-1}{\big|}\}.
\end{align*}

\medskip
\noindent
{\it Step 4.} \textit{Derivation of $p_{k,1}$} \\
According to Cournat's Theorem~\ref{thm:courant}, we take a Jordan curve $\Gamma_1$ such that $\Gamma_1 \subset \C \setminus \overline{\inte \Gamma_0}$ and we choose a continuous mapping $ \Phi_1 : \overline{\inte \Gamma_1} \to \D$ which is conformal on $\inte \Gamma_1 $   satisfying $ \Phi_1 (0)=0$ and
\begin{align}\label{eq:conf-cond_p-1}
 | \Phi_0 (z)- \Phi_1 (z)| < \frac{\delta_k}{2}
 \qquad \text{ for all }  z \in \overline{\inte \Gamma_0}.
\end{align}
The construction of $\Gamma_1$ and hence of $\Phi_1$ depends on the properties of $a_0$. We note that numerical methods for conformal mappings is a research area itself and, thus, it is beyond the scope of this paper. We refer to the comments provided at the end of Section~\ref{sec:approx}.

To derive the polynomial $p_{k,1}$ we apply Runge's Little Theorem~\ref{thm:runge} to the holomorphic function $\Phi_1 : \inte \Gamma_1 \to \D$  so that
\begin{align*}
 |p_{k,1}(z) - \Phi_1(z)|< \frac{\delta_k}{2}
 \qquad \text{ for all } z \in \overline{\inte \Gamma_0}.
\end{align*}
Then, for all $z \in \overline{\inte \Gamma_0}${,} one has
\begin{align*}
 |p_{k,1}(z) - \Phi_0(z)| \leq  |p_{k,1}(z) - \Phi_1(z)| + | \Phi_1(z) - \Phi_0(z)|< \delta_k .
\end{align*}
Further, as shown in the proof of Theorem~\ref{thm:walsh_jordan_curve_poly} (cf.\eqref{eq:walsh-jordan-kniff}),  for all $z\in \Gamma_0$ it follows
\begin{align*}
\left|\sum_{l=0}^{n_k-1} c_l \, \big(\Phi_0(z)^l- p_{k,1}(z)^l\big)        \right|
<\frac{\e}{4  \|T\|_{M,\infty}}.
\end{align*}

\medskip
\noindent
{\it Step 5.} \textit{Derivation of $p_{k,2}$}\\
We consider the mapping $\Phi_0:\overline{\inte \Gamma_0} \to \overline{\D}$ and the residue of $\tfrac{1}{\Phi_0}$ at zero, i.e.{,}
$$a_{-1} = \frac{1}{2\pi i} \int_{\partial B_r(0)} \frac{1}{\Phi_0(\xi)} \di{\xi}, \qquad r >0 \, \, \text{such that } \overline{B_r(0)} \subset \inte \Gamma_0. $$ 
Moreover, let
$\Psi_0\colon \inte \Gamma_0\to \C$,
\begin{align}\label{eq:def-phi_e-walsh-curve-0}
\Psi_0(z) = \begin{cases} 
\frac{1}{\Phi_0(z)} -  \frac{a_{-1}}{z}  & z \in \inte \Gamma_0\setminus \{0\}\\
\frac{1}{2\pi i} \int_{\partial B_r(0)} \frac{1}{\xi \Phi_0(\xi)} \di{\xi} & z=0,
\end{cases} 
\end{align}
denote the holomorphic extension, which is continuous on $\overline{ \inte \Gamma_0}$. To determine $p_{k,2}${,} we take another Jordan curve $\Gamma_2$ such that
$\Gamma_2 \subset \C \setminus \overline{\inte \Gamma_0}$ and we choose a continuous mapping $ \Phi_2 : \overline{\inte \Gamma_2} \to \D$ which is conformal on $\inte \Gamma_2 $   satisfying $ \Phi_2 (0)=0$ and
\begin{align}\label{eq:conf-cond_p-2}
 | \Phi^{-1}_0 (w) - \Phi^{-1}_2 (w)| < \frac{\delta_k}{2L_{\Psi_0}}
 \qquad \text{ for all }  w \in \overline{\D},
\end{align}
where $L_{\Psi_0}>0$ denotes the Lipschitz constant of $\Psi_0$\footnote{By Kellogg's Theorem, as $a_0$ and its extension $\Gamma_0$ are twice continuously differentiable, the function $\Psi_0$ satisfies a Lipschitz condition on $\Gamma_0$.}. 
Then, for all $  z =\Phi_2^{-1}(w) \in \overline{\inte \Gamma_0}${,} we have
\begin{align*}
 | \Psi_0 (z)- \Psi_0\circ\Phi^{-1}_0\circ\Phi_2 (z)| \leq L_{\Psi_0} \,  | \Phi^{-1}_2 (w) - \Phi^{-1}_0 (w)| < \frac{\delta_k}{2}.
\end{align*}
Next, we apply Runge's Little Theorem to the holomorphic function
$$ \Psi_0\circ\Phi^{-1}_0\circ\Phi_2 :\inte \Gamma_2  \to \C $$
and determine the polynomial $p_{k,2}$ such that
\begin{align*}
| p_{k,2} (z) - \Psi_0\circ\Phi^{-1}_0\circ\Phi_2 (z)| < \frac{\delta_k}{2}
 \qquad \text{ for all }  z \in \overline{\inte \Gamma_0},
\end{align*}
Thus, as shown in the proof of Theorem~\ref{thm:walsh_jordan_curve_poly} (similar to \eqref{eq:walsh-jordan-kniff}), for all $z\in \Gamma_0$ it follows
\begin{align*}
\left|\sum_{l=1}^{n_k-1} c_{-l} \, 
\Big( \, \big( \tfrac{a_{-1}}{z} +p_{k,2}(z)\big)^l -  \big(\tfrac{1}{\Phi_0(z)}\big)^l \Big)
\right|
<\frac{\e}{4  \|T\|_{M,\infty}}.
\end{align*}

\medskip
\noindent
{\it Step 6.} \textit{Derivation of $p_{k,3}$}\\
In this step we will apply Runge's Little Theorem to the holomorphic function $z \mapsto \tfrac{a_{-1}}{z}$ on $a_0(\p)$. To this end, let 
\begin{align*}
 \nu:= \max_{z \in a_0(\p)} \left|\tfrac{1}{z} \right|, \qquad  
 \mu_k:= \max_{z \in a_0(\p)} \left|p_{k,2}(z) \right|.
\end{align*}
Then, we take $\eta_k>0$ such that
\begin{align*}
\eta_k  \cdot \frac{1-(2\eta_k + |a_{-1}|\nu +\mu_k)^{n_k-1}}{1-(2\eta_k + |a_{-1}|\nu +\mu_k)}   < \frac{\e}{4 (n_k-1)c_{\max} \|T\|_{M,\infty}}.
\end{align*}
The polynomial $p_{k,3}$ is then obtained by applying Runge's Little Theorem  such that for all $z \in a_0(\p)$ one has
\begin{align*}
|\tfrac{a_{-1}}{z} -p_{k,3}(z)|< \eta_k . 
\end{align*}
Then, as shown in the proof of Theorem~\ref{thm:walsh-arc}, we get
\begin{align*}
 \left|\sum_{l=1}^{n_k-1} c_{-l} \, \Big( \,    \big(p_{k,2}(z) + p_{k,3}(z) \big)^l
- \big( \tfrac{a_{-1}}{z} +p_{k,2}(z)\big)^l  \Big)
\right| < \frac{\e}{4   \|T\|_{M,\infty}}.
\end{align*}
Thus, we have
\begin{align*}
| (T^{-1}f)_k\circ a_0^{-1}(z) -  p_k(z)| < \frac{\e}{  \|T\|_{M,\infty}}.
\end{align*}

\medskip
\noindent
{\it Step 7.} \textit{Monomial representation of $p$}\\
Finally,  we determine the monomial representation of the polynomial 
\begin{align}\label{def:poly_walsh_S1}
p(z) = \sum_{k=1}^n p_k(z^n+a_{n-1}z^{n-1}+ \cdots + a_1\,z) z^{k-1}   =:\sum_{l=0}^N \beta_l z^l.
\end{align}
This construction yields the following

\begin{theorem}
{Let $\p$ be a Jordan arc having a twice continuously differentiable parametrization.} Assume that $(A,b)\in C^1_{n,n}(\p)\times C_n^1(\p)$ satisfies (N1), (N2) and (S1). Let $\varepsilon>0$,  $f \in \operatorname{C}_n(\p)$ and $a_0 \in C^2(\p)$ be a Jordan arc in the complex plane such that $(T^{-1}f)_k \circ a_0^{-1} \in \operatorname{Lip}(a_0(\p))$ for all $k=1,...,n$. Then, for
\begin{align*}
T=N+1 \qquad \text{and} \qquad  u=(u_0,...,u_{T-1}) = (\beta_N,...,\beta_{0})
\end{align*}
defined in \eqref{def:poly_walsh_S1}{,} one has
\begin{align*}
\| \varphi(T,u,0) - f\|_{\infty}  < \varepsilon.
\end{align*}
In particular, the result is true if $(A,b) \in C^1_{n,n}(\p)\times C^1_n(\p)$ and $f\in \operatorname{Lip}_n(\p)$.
\end{theorem}

\subsection{Methods for condition~(S2)}\label{sec:meth_B}

In this section we consider methods to construct suitable inputs based on the necessary conditions (N1), (N2) and the sufficiency condition~(S2). Let $f \in \operatorname{C}_n(\p)$  and $\varepsilon>0$ be given. After applying a change of coordinates $T(\theta)${,} we consider the pair
\begin{align*}
\tilde A(\theta):= T(\theta)^{-1} A(\theta)T(\theta)=\begin{pmatrix}
 \lambda_1(\theta) & & \\
 & \ddots &\\
 & & \lambda_n(\theta)
\end{pmatrix},\quad
\tilde b(\theta):=T(\theta)^{-1}b(\theta) = \begin{pmatrix}
 1\\
\vdots\\
1
\end{pmatrix},
\end{align*}
where $\lambda_1,...,\lambda_n$ denote the distinct eigenvalue Jordan arcs. Writing the solution formula in terms of a polynomial $p${,} we get
\begin{align*}
\| \varphi(T,u,0) - f\|_{\infty}  &= \| p(A)b - f\|_{\infty} = \|T ( p(\tilde A)\tilde b - T^{-1}f)\|_{\infty} \\&\leq \|T\|_{M,\infty} \,\| p(\tilde A)\tilde b - T^{-1}f\|_{\infty}\\
&= \|T\|_{M,\infty} \,\left\| \begin{pmatrix} p(\lambda_1) - (T^{-1}f)_1\\\vdots\\
                        p(\lambda_n) - (T^{-1}f)_n  
                         \end{pmatrix}
\right\|_{\infty},
\end{align*}
where $(T^{-1}f)_k$ denotes the $k$th component of the vector-valued function $\theta \mapsto T^{-1}(\theta)f(\theta)\in \C^n$.
Thus, we conclude that
\begin{align*}
\| \varphi(T,u,0) - f\|_{\infty}  <  \e {,}
\end{align*}
if we can  find a polynomial  $p$ so that for all $k=1,...,n$ it holds
\begin{align*}
\sup_{z \in \lambda_k(\p)} | p(z) - (T^{-1}f)_k \circ \lambda_k^{-1}(z)|  < \frac{\e}{\|T\|_{M,\infty}}.
\end{align*}
Compared to the methods for (S1), we are a looking for a single polynomial working for every component of the target function. With other words, the polynomial $p$ solves $n$ approximation problems simultaneously, where each problem is defined on an individual domain. Inspired by \cite{andrievskii2005polynomial}, using the ansatz
\begin{equation*}
p(z) = \sum_{k=1}^n p_k(z)q_k(z),
\end{equation*}
we will split the problem into $n$ subproblems. Postponing the technicalities for a moment, on {the} one hand we will determine $p_1,...,p_n$ so that for each $k\in \{1,...,n\}$ the polynomial $p_k$ approximates the function $(T^{-1}f)_k \circ \lambda_k^{-1}$ over the image of $\lambda_k\colon \p \to \C$. Let $\lambda_k(\p)$ denote the image of $\lambda_k$.  On {the} other hand, the polynomials $q_1,...,q_n$ {are} derived by approximating the indicator functions
\begin{align*}
h_k\colon \lambda(\p) \to \C , \quad   h_k(z)= 
 \begin{cases}
  1 & \text{ if } z \in \lambda_k (\p) \\
  0 & \text{ if } z \in \lambda(\p) \setminus \lambda_k(\p),
 \end{cases}
\end{align*}
where 
 $$\lambda(\p) = \bigcup_{k=1}^n \lambda_k(\p). $$
Note that the assumptions imply that the sets $\lambda_1(\p),..., \lambda_n(\p)$ are pairwise disjoint and so  $h_1,...,h_n$ are holomorphic functions.

A precise outline for the construction of the polynomials $p_1,...,p_n$ and $q_1,...,q_n$ is presented in the following. In doing so, we will distinguish two cases. Exemplary we will discuss the construction process for polynomials $p_k$ and $q_k$, for some $k \in \{1,...,n\}$.

\subsubsection*{Case~1: $\lambda_k(\p)=[a,b]\subset \R$ defines a compact interval}

We suppose the function $z \mapsto (T^{-1}f)_k\circ \lambda_k^{-1}(z)$ is Lipschitz continuous\footnote{This is the case, e.g. if the pair $(A,b)$ is in $C^1_{n,1}(\p)\times C_n^1(\p)$ {and $f \in \operatorname{Lip}_n(\p)$}.}. Let $L_k{>}0$ denote the Lipschitz constant. Then, by Theorem~\ref{thm:gzyl} the condition
\begin{align}\label{eq:S2_p-k-Bernstein}
\sup_{z \in \lambda_k(\p)} | p_k(z) - (T^{-1}f)_k \circ \lambda_k^{-1}(z)|  < \frac{\e}{\|T\|_{M,\infty}}
\end{align}
is achieved by taking $p_k= B_{(T^{-1}f)_k\circ \lambda_k^{-1},n_k}$ to be the Bernstein polynomial of the function $z \mapsto (T^{-1}f)_k\circ  \lambda_k^{-1}(z)$ of degree $n_k \geq 3$ satisfying
\begin{equation*}
\sqrt{2} \left( 4\,M_{(T^{-1}f)_k\circ \lambda_k^{-1}} + \frac{(b-a)\,L_{k}}{2} \right) {\sqrt\frac{\log n_k}{n_k}} < \frac{\e}{\|T\|_{M,\infty}}.
\end{equation*}
After we determined $p_1,...,p_n${,} we define
$$\alpha_{k,l}:= \sup_{\theta \in \p}|B_{(T^{-1}f)_l\circ \lambda_l^{-1},n_l}(\lambda_k(\theta))|.$$
The polynomial $q_k$ is derived by applying Runge's Little Theorem to the holomorphic function $h_k$ such that for each $k=1,...,n$ it holds
\begin{align}\label{eq:S2_q-k-Bernstein}
\sup_{z \in \lambda_k(\p)} | q_k(z) - h_k(z))| <  \frac{\varepsilon}{3\|T\|_{M,\infty} \, \sum_{l=1}^n \alpha_{k,l}}.
\end{align}
{For every $k=1,...,n$, based on the conditions~\eqref{eq:S2_p-k-Bernstein} and \eqref{eq:S2_q-k-Bernstein},  we get}
\begin{align*}
| (T^{-1}&f)_k \circ \lambda_k^{-1}(z) -p(z)| \leq | (T^{-1}f)_k \circ \lambda_k^{-1}(z)- p_k(z)q_k(z)| + \sum_{l\neq k} |q_l(z)|\, |p_l(z)|\\
&\leq  | (T^{-1}f)_k \circ \lambda_k^{-1}(z)- p_k(z)| + |p_k(z)||1-q_k(z)| + \sum_{l\neq k} |q_l(z)|\, |p_l(z)|\\
&\leq  \frac{\e}{3 \|T\|_{M,\infty}} +  \frac{\e \, \alpha_{k,k}}{3 \|T\|_{M,\infty}\sum_{l=1}^n \alpha_{k,l}} + \sum_{l \neq k} \frac{\e \, \alpha_{k,l}}{3 \|T\|_{M,\infty} \sum_{l=1}^n \alpha_{k,l}} \leq  \frac{\e}{ \|T\|_{M,\infty}}.
\end{align*}
Consequently, by considering the monomial representation
\begin{align}\label{def:poly_Bernstein_S2}
p(z)= \sum_{k=1}^n q_k(z)\,  B_{(T^{-1}f)_k\circ \lambda_k^{-1},n_k} (z)  =\sum_{k=0}^N \beta_k z^k
\end{align}
we get  the following {result.}

\begin{theorem}
Let $(A,b) \in C_{n,n}(\p)\times C_n(\p)$ satisfy (N1), (N2) and (S2). Let $\varepsilon>0$ and suppose that $f \in C_n(\p)$ is such that $(T^{-1}f)_k \circ \lambda_k^{-1} \in \operatorname{Lip}(\lambda_k(\p))$ for all $k=1,...,n$ and assume that $\lambda_1(\p),...,\lambda_n(\p)$ are compact real intervals.  Then, for
\begin{align*}
T=N+1 \qquad \text{and} \qquad  u=(u_0,...,u_{T-1}) = (\beta_N,...,\beta_{0})
\end{align*}
defined in \eqref{def:poly_Bernstein_S2}{,} one has
\begin{align*}
\| \varphi(T,u,0) - f\|_{\infty} < \varepsilon.
\end{align*}
In particular, the result is true if $(A,b) \in C^1_{n,n}(\p)\times C^1_n(\p)$ and $f\in \operatorname{Lip}_n(\p)$.
\end{theorem}

Next we consider the {more} general case where the eigenvalue functions define   Jordan arcs.

\subsubsection*{Case~2: $\lambda_1(\p),...,\lambda_n(\p)$ define Jordan arcs}

The construction process is similar to the methods for (S1). Hence, in order to avoid a lengthy presentation we will go through the analogous construction steps more quickly.

\medskip
\noindent
{\it Step 1.} \textit{Extension of the Jordan arcs and the target function}\\
Without loss of generality we will assume that none of the Jordan arcs contains the origin.
In a first step we extent the Jordan arcs $\lambda_1(\p), ...,\lambda_n(\p)$ to Jordan curves $\Gamma_1,..., \Gamma_n$ such that the origin lies inside {of} each Jordan curve. Then {we} take a continuous bijective mappings $\Phi_k : \overline{\inte \Gamma_k} \to \overline{\mathbb D}$ which are conformal on $\inte \Gamma_k$ so that $\Phi_k(0)=0$. The existence of such mappings is ensured by Caratheodory's Theorem, cf. \cite[Thm.~13.2.3]{greene_krantz_1997}.

Analogous to Step~1 in the methods for (S1), we extend the functions  $ z \mapsto (T^{-1}f)_k\circ \lambda_k^{-1}(z)$ to $\Gamma_k$, where $T\in C_{n,n}(\p)$ denotes the changes of coordinates introduced at the beginning of Section~\ref{sec:meth_B}. As a result{,} we obtain continuous functions
\[
h_k := (T^{-1}f)_k\circ \lambda_k^{-1} \circ \Phi_k^{-1}  \colon \partial\mathbb{D} \to \mathbb C{, \quad k=1,...,n}.
\]
If  $\lambda_1,...,\lambda_n \in C^2(\p)$ and $\p$ has a parametrization which {is} twice continuously differentiable and the extensions $\Gamma_1,...,\Gamma_n$ have the same properties, then by Kellogg's Theorem \cite[Theorem~4.3]{garnett_marshall_2005}, the functions $h_k$ satisfy a Lipschitz condition provided {the} $(T^{-1}f)_k$ {are} Lipschitz continuous\footnote{This is the case if $(A,b) \in C_{n,n}^1(\p) \times C_n^1(\p)$ {and $f\in \operatorname{Lip}_n(\p)$}.}. Let $L_k>0$ denote the corresponding Lipschitz constant.

\medskip
\noindent
{\it Step 2.} \textit{Ansatz for the polynomial $p$}\\
As in the previous case, to derive a suitable polynomial, we set
\begin{align*}
p(z)= \sum_{k=1}^n p_k(z)\, q_{k}(z).
\end{align*}
Based on the exploration in Section~\ref{sec:appox_walsh} (in particular~\eqref{eq:def-poly-walsh-arc}), for each $k =1,...,n$ we  make the ansatz
\begin{align*}
 p_k(z) = \sum_{k=0}^{n_k-1} c_l \, \big(p_{k,1}(z)\big)^l + \sum_{l=1}^{n_k-1} c_{-l} \,   \big(p_{k,2}(z) + p_{k,3}(z)\big)^l,
\end{align*}
where  
\begin{align*}
 c_l:= \frac{n_k - |l|}{n_k}\,  \widehat{h_k} (l),\qquad \, l = -n_k+1,....,n_k-1
\end{align*}
denote the Fourier coefficients corresponding to the function $h_k$. Again we take a simple but instructive step
\begin{align*}
| (T^{-1}f)_k\circ \lambda_k^{-1}(z) &-  p_k(z)| \\
&\leq 
\left|   (T^{-1}f)_k\circ a_0^{-1}  (z) -  \sum_{l=0}^{n_k-1} c_l \, \big(\Phi_k(z)\big)^l  - 
\sum_{l=1}^{n_k-1} c_{-l} \,   \big(\tfrac{1}{\Phi_k(z)}\big)^l\right|\\
&\quad+ \left|\sum_{l=0}^{n_k-1} c_l \, \big(\Phi_k(z)^l- p_{k,1}(z)^l\big)        \right|
\\ &\quad
+ \left|\sum_{l=1}^{n_k-1} c_{-l} \, \Big( \, \big( \tfrac{a_{-1,k}}{z} +p_{k,2}(z)\big)^l -  \big(\tfrac{1}{\Phi_k(z)}\big)^l \Big)
\right|
\\ &\quad
+ \left|\sum_{l=1}^{n_k-1} c_{-l} \, \Big( \,    \big(p_{k,2}(z) + p_{k,3}(z) \big)^l
- \big( \tfrac{a_{-1,k}}{z} +p_{k,2}(z)\big)^l  \Big)
\right|,
\end{align*}
where
$$a_{-1,k} = \frac{1}{2\pi i} \int_{\partial B_r(0)} \frac{1}{\Phi_k(\xi)} \di{\xi}, \qquad r >0 \, \, \text{such that } \overline{B_r(0)} \subset \inte \Gamma_0$$ 
denotes the residue of $\tfrac{1}{\Phi_k}$ at zero.  

\medskip
\noindent
{\it Step 3.} \textit{Derivation of $n_k \in \N$}\\
Applying Weierstrass Second Theorem~\ref{thm:weier-circ}~(b) to the function
\[
 h_k = (T^{-1}f)_k\circ \lambda_k^{-1} \circ \Phi_k^{-1}  \colon \partial\mathbb{D} \to \mathbb C
\]
we get
\begin{align*}
\sup_{w \in \partial \D}\left|h_k (w) -  \sum_{l=0}^{n_k-1} c_l \, w^l  - 
\sum_{l=1}^{n_k-1} c_{-l} \,   \big(\tfrac{1}{w}\big)^l\right| 
< 2\sqrt{2}\pi L_{k} \frac{\ln n_k }{n_k} .
\end{align*}
As $\Phi_k$ is one-to-one and onto and taking $n_k \in \N$ so that 
\begin{align*}
  2\sqrt{2}\pi L_{k} \frac{\ln n_k }{n_k}   \leq  \frac{\e}{12\,\|T\|_{M,\infty}},
\end{align*} 
we have
\begin{align*}
\sup_{z \in \Gamma_k} \left| (T^{-1}f)_k\circ \lambda_k^{-1} (z) -  \sum_{l=0}^{n_k-1} c_l \, \big(\Phi_k(z)\big)^l  - 
\sum_{l=1}^{n_k-1} c_{-l} \,   \big(\tfrac{1}{\Phi_k(z)}\big)^l\right| 
< \frac{\e}{12 \, \|T\|_{M,\infty}}. 
\end{align*}

Hence, it remains to determine the polynomials $p_{k,1},p_{k,2},p_{k,3}$.
This will be carried out in the following three steps. To this end, we pick $\delta_{k}>0$ so that
\begin{align*}
  \sum_{l=0}^{n_k-1} \left( (1+\delta_k)^l -1 \right) < \frac{\e}{12\,c_{\max} \| T\|_{M,\infty}}, \qquad c_{\max} := \max \{{|}c_{-n_k+1}{|},...,{|}c_{n_k-1}{|}\}.
\end{align*}

\medskip
\noindent
{\it Step 4.} \textit{Derivation of $p_{k,1}$} \\
According to Cournat's Theorem~\ref{thm:courant}, we take a Jordan curve $\Gamma_{k,1}$ such that $ \Gamma_{k,1} \subset \C \setminus \overline{\inte \Gamma_k}$ and we choose a continuous mapping $ \Phi_{k,1} : \overline{\inte \Gamma_{k,1}} \to \overline{\D}$ which is conformal on $\inte \Gamma_{k,1}$ satisfying $ \Phi_{k,1} (0)=0$ and
\begin{align}\label{eq:conf-cond_p,k-1}
 | \Phi_k (z)- \Phi_{k,1} (z)| < \frac{\delta_k}{2}
 \qquad \text{ for all }  z \in \overline{\inte \Gamma_k}.
\end{align}
The construction of $\Gamma_{k,1}$ and hence of $\Phi_{k,1}$ depends on the properties of $\lambda_k$. Recall that constructive methods for conformal mappings are beyond the scope of this paper. We refer to the comments provided at the end of Section~\ref{sec:approx}. 

To derive the polynomial $p_{k,1}${,} we apply Runge's Little Theorem~\ref{thm:runge} to the conformal mapping $\Phi_{k,1}$  so that
\begin{align*}
 |p_{k,1}(z) - \Phi_{k,1}(z)|< \frac{\delta_k}{2}
 \qquad \text{ for all } z \in \overline{\inte \Gamma_k}.
\end{align*}
Then, for all $z \in \overline{\inte \Gamma_k}$ one has
\begin{align*}
 |p_{k,1}(z) - \Phi_k(z)| \leq  |p_{k,1}(z) - \Phi_{k,1}(z)| + | \Phi_{k,1}(z) - \Phi_k(z)|< \delta_k .
\end{align*}
Further, as shown in the proof of Theorem~\ref{thm:walsh_jordan_curve_poly} (cf.\eqref{eq:walsh-jordan-kniff}),  for all $z\in \Gamma_k$ it follows
\begin{align*}
\left|\sum_{l=0}^{n_k-1} c_l \, \big(\Phi_k(z)^l- p_{k,1}(z)^l\big)        \right|
<\frac{\e}{12  \|T\|_{M,\infty}}.
\end{align*}

\medskip
\noindent
{\it Step 5.} \textit{Derivation of $p_{k,2}$}\\
We consider the holomorphic function
$\Psi_k\colon \inte \Gamma_k\to \C$,
\begin{align}\label{eq:def-phi_e-walsh-curve-k}
\Psi_k(z) = \begin{cases} 
\frac{1}{\Phi_k(z)} -  \frac{a_{-1,k}}{z}  & z \in \inte \Gamma_k\setminus \{0\}\\
\frac{1}{2\pi i} \int_{\partial B_r(0)} \frac{1}{\xi \Phi_k(\xi)} \di{\xi} & z=0,
\end{cases} 
\end{align}
which is continuous on $\overline{ \inte \Gamma_k}$. To determine $p_{k,2}${,} we take another Jordan curve $\Gamma_{k,2}$ such that $\Gamma_{k,2} \subset \C \setminus \overline{\inte \Gamma_k}$ and we choose a continuous mapping $ \Phi_{k,2} : \overline{\inte \Gamma_{k,2}} \to \overline{\D}$ which is conformal on $\inte \Gamma_{k,2}$ satisfying $ \Phi_{k,2} (0)=0$ and
\begin{align}\label{eq:conf-cond_p,k-2}
 | \Phi^{-1}_k (w) - \Phi^{-1}_{k,2} (w)| < \frac{\delta_k}{2L_{\Psi_k}}
 \qquad \text{ for all }  w \in \overline{\D},
\end{align}
where $L_{\Psi_k}>0$ denotes the Lipschitz constant of $\Psi_k$. 
Then, for all $  z =\Phi_{k,2}^{-1}(w) \in \overline{\inte \Gamma_k}$ we have
\begin{align*}
 | \Psi_k (z)- \Psi_k\circ\Phi^{-1}_0\circ\Phi_{k,2} (z)| \leq L_{\Psi} \,  | \Phi^{-1}_{k,2} (w) - \Phi^{-1}_k (w)| < \frac{\delta_k}{2}.
\end{align*}
Next, we apply Runge's Little Theorem to the holomorphic function
$$ \Psi_k\circ\Phi^{-1}_k\circ\Phi_{k,2} :\inte \Gamma_{k,2}  \to \C $$
and determine the polynomial $p_{k,2}$ such that
\begin{align*}
| p_{k,2} (z) - \Psi_k\circ\Phi^{-1}_k\circ\Phi_{k,2} (z)| < \frac{\delta_k}{2}
 \qquad \text{ for all }  z \in \overline{\inte \Gamma_k},
\end{align*}
Thus, as shown in the proof of Theorem~\ref{thm:walsh_jordan_curve_poly} (similar to \eqref{eq:walsh-jordan-kniff}), for all $z\in \Gamma_k$ it follows
\begin{align*}
\left|\sum_{l=1}^{n_k-1} c_{-l} \, 
\Big( \, \big( \tfrac{a_{-1,k}}{z} +p_{k,2}(z)\big)^l -  \big(\tfrac{1}{\Phi_k(z)}\big)^l \Big)
\right|
<\frac{\e}{12  \|T\|_{M,\infty}}.
\end{align*}

\medskip
\noindent
{\it Step 6.} \textit{Derivation of $p_{k,3}$}\\
In this step we will apply Runge's Little Theorem to the holomorphic function $z \mapsto \tfrac{a_{-1,k}}{z}$ on $\lambda_k(\p)$. More precisely, let 
\begin{align*}
 \nu_k:= \max_{z \in \lambda_k(\p)} \left|\tfrac{1}{z} \right|, \qquad  
 \mu_k:= \max_{z \in \lambda_k(\p)} \left|p_{k,2}(z) \right|.
\end{align*}
Then, we take $\eta_k>0$ such that
\begin{align*}
\eta_k  \cdot \frac{1-(2\eta_k + |a_{-1,k}|\nu_k +\mu_k)^{n_k-1}}{1-(2\eta_k + |a_{-1,k}|\nu_k +\mu_k)}   < \frac{\e}{12 (n_k-1)c_{\max} \|T\|_{M,\infty}}.
\end{align*}
The polynomial $p_{k,3}$ is then obtained by applying Runge's Little Theorem  such that for all $z \in \lambda_k(\p)$ one has
\begin{align*}
|\tfrac{a_{-1,k}}{z} -p_{k,3}(z)|< \eta_k . 
\end{align*}
Then, as shown in the proof of Theorem~\ref{thm:walsh-arc}, we get
\begin{align*}
 \left|\sum_{l=1}^{n_k-1} c_{-l} \, \Big( \,    \big(p_{k,2}(z) + p_{k,3}(z) \big)^l
- \big( \tfrac{a_{-1}}{z} +p_{k,2}(z)\big)^l  \Big)
\right|<\frac{\e}{12  \|T\|_{M,\infty}}.
\end{align*}
Thus, we have
\begin{align}\label{eq:S2-target-p_k-jordan}
| (T^{-1}f)_k\circ \lambda_k^{-1}(z) &-  p_k(z)| < \frac{\e}{3 \|T\|_{M,\infty}}.
\end{align}

\medskip
\noindent
{\it Step 7.} \textit{Derivation of $q_{k}$}\\
Let $\lambda_k(\p)$ denote the image of $\lambda_k\colon \p \to \C$ and 
 $$\lambda(\p) = \bigcup_{k=1}^n \lambda_k(\p). $$
The indicator function on $\lambda (\p)$ corresponding to the Jordan arc $\lambda_k$ is denoted by
\begin{align*}
h_k\colon \lambda(\p) \to \C , \quad   h_k(z)= 
 \begin{cases}
  1 & \text{ if } z \in \lambda_k (\p) \\
  0 & \text{ if } z \in \lambda(\p) \setminus \lambda_k(\p).
 \end{cases}
\end{align*}
Recall that the assumptions imply that $\lambda_1(\p),..., \lambda_n(\p)$ are pairwise disjoint and so  $h_k$ is a holomorphic function. Next, we define
$$\alpha_{k,l}:= \sup_{\theta \in \p}|p_l(\lambda_k(\theta))|.$$
The polynomial $q_k$ is then derived by applying Runge's Little Theorem to  $h_k$ such that
\begin{align}\label{eq:S2_q-k-jordan}
\sup_{z \in \lambda_k(\p)} | q_k(z) - h_k(z))| <  \frac{\varepsilon}{3\|T\|_{M,\infty} \, \sum_{l=1}^n \alpha_{k,l}}.
\end{align}
Based on the conditions~\eqref{eq:S2-target-p_k-jordan} and \eqref{eq:S2_q-k-jordan}, for every $k=1,...,n$ we get
\begin{align*}
| (T^{-1}&f)_k \circ \lambda_k^{-1}(z) -p(z)| \leq | (T^{-1}f)_k \circ \lambda_k^{-1}(z)- p_k(z)q_k(z)| + \sum_{l\neq k} |q_l(z)|\, |p_l(z)|\\
&\leq  | (T^{-1}f)_k \circ \lambda_k^{-1}(z)- p_k(z)| + |p_k(z)||1-q_k(z)| + \sum_{l\neq k} |q_l(z)|\, |p_l(z)|\\
&\leq  \frac{\e}{3 \|T\|_{M,\infty}} +  \frac{\e \, \alpha_{k,k}}{3 \|T\|_{M,\infty}\sum_{l=1}^n \alpha_{k,l}} + \sum_{l \neq k} \frac{\e \, \alpha_{k,l}}{3 \|T\|_{M,\infty} \sum_{l=1}^n \alpha_{k,l}} \leq  \frac{\e}{ \|T\|_{M,\infty}}.
\end{align*}

\medskip
\noindent
{\it Step 8.} \textit{Monomial representation of $p$}\\
Finally, using the polynomials  $p_1,...,p_n$ and $q_1,...,q_n${,} we determine the monomial representation of the polynomial
\begin{align}\label{def:poly_walsh_S2}
p(z) = \sum_{k=1}^n p_k(z) q_k(z)   =:\sum_{l=0}^N \beta_l z^l.
\end{align}
This construction yields the following {result.}

\begin{theorem}
{Let $\p$ be a Jordan arc having a twice continuously differentiable parametrization.} Assume that $(A,b)\in C^1_{n,n}(\p)\times C_n^1(\p)$ satisfies (N1), (N2) and (S1). Let $\varepsilon>0$,  $f \in \operatorname{C}_n(\p)$ and $\lambda_1,...,\lambda_n \in C^2(\p)$ be  Jordan arcs in the complex plane such that $(T^{-1}f)_k \circ \lambda_k^{-1} \in \operatorname{Lip}(\lambda_k(\p))$ for all $k=1,...,n$. Then, for
\begin{align*}
T=N+1 \qquad \text{and} \qquad  u=(u_0,...,u_{T-1}) = (\beta_N,...,\beta_{0})
\end{align*}
defined in \eqref{def:poly_walsh_S2}{,} one has
\begin{align*}
\| \varphi(T,u,0) - f\|_{\infty}  < \varepsilon.
\end{align*}
In particular, the result is true if $(A,b) \in C^1_{n,n}(\p)\times C^1_n(\p)$ and $f\in \operatorname{Lip}_n(\p)$.
\end{theorem}

\section{Continuous-time single-input systems}\label{sec:poly_approx_cont}

In this section{,} we investigate how to get constructive methods for continuous-time single-input systems. First, we note that a direct application of the discrete-time methods above is not immediately possible. The reason is that, for $u\in L^1([0,T])$  the solution to
\begin{align*}
\frac{\partial x}{\partial t}
(t,\theta)=A(\theta)x(t,\theta)+b(\theta)u(t), \quad x(0,\theta)=x_0(\theta)
\end{align*}
is given by
\begin{align*}
\varphi(T,u,x_0) (\theta) = {\rm e}^{TA(\theta)} x_0(\theta) +  \int_0^T {\rm e}^{(T-s)A(\theta)}b(\theta)u(s)\di{s},
\end{align*}
which is not of the form $p(A(\theta))b(\theta)$ for some polynomial $p$.  In the following{,} we present an approach that enables us to trace the input construction problem for continuous-time systems back to the methods just presented in the previous section. The basic idea simply is to approximate the integral in the solution formula by a polynomial using piecewise constant inputs. The analysis will be carried out in six steps. Before we start, there are two comments in order. For continuous-time systems, the final time can be chosen arbitrarily, i.e. we can fix some $T>0$ in advance. As a consequence, for continuous-time systems we are not restricted to start at the origin. Indeed, let $x_0 \in C_n(\p)$ denote the initial states and let $f\in C_n(\p)$ denote the final states, then using $\tilde f \in C_n(\p)$, defined by
\[
 \tilde f(\theta) := f(\theta) -  {\rm e}^{TA(\theta)} x_0(\theta) 
\]
an input $u \in L^1([0,T])$ satisfies $\|\varphi(T,u,x_0)-f\|_{\infty}<\e$ if and only if $\|\varphi(T,u,0)-\tilde f\|_{\infty}<\e$. 

Suppose the pair $(A,b)$ satisfies the conditions (N1), (N2) and (S2) and {let} $f \in C_n(\p)$ and $\e>0$ {be} given. As $T>0$ can chosen arbitrarily, we set  $T = N\tau$ for some $\tau>0$ and some $N \in \N$ that will be specified in the construction process.

\medskip
\noindent
{\it Step~1: Diagonalize $A(\theta)$} \\
By assumptions (N2) and (S2) there is a continuous change of coordinates $T(\theta)$ such that
\begin{equation*}
T(\theta)^{-1} A(\theta) T(\theta) = 
\begin{pmatrix}
\lambda_1(\theta) & & 0\\
  & \ddots & \\
0 & & \lambda_n(\theta)
\end{pmatrix}
\quad 
\text{ and } \quad 
T(\theta)^{-1} b(\theta)= 
\begin{pmatrix}
1\\
 \vdots  \\
1
\end{pmatrix}.
\end{equation*}
The solution can be written as
\begin{align*}
\varphi(N\tau,u,0) (\theta) 
= T(\theta) 
\begin{pmatrix}
\tilde \varphi_1(N\tau,u,0)(\theta) \\
\vdots\\
\tilde \varphi_n(N\tau,u,0)(\theta)
\end{pmatrix}, 
\end{align*}
where
\begin{align*}
\tilde \varphi_k(N\tau,u,0)(\theta)= \int_0^{N\tau} {\rm e}^{(N\tau-s)\lambda_k(\theta)}    u(s) \di{s} , \quad k=1,...,n.
\end{align*}

\medskip
\noindent
{\it Step~2: Take piecewise constant inputs}\\
Let $\tau>0$ (not yet specified) and consider the following partition
$$[0,N\tau] = [0,\tau) \cup [\tau, 2\tau) \cup \cdots \cup [(N-1)\tau, N\tau],$$
where every interval $I_l:=[l\tau, (l+1)\tau)$ has length $\tau$. Further, we consider inputs $u\colon [0,T] \to \C$ that are constant on every interval, i.e.
\begin{align}\label{def:u}
u\vert_{I_l}(t):= u_l \in \C, \quad l =0,...,N-1.
\end{align} 
Let $\mathbf{1}_{I_l}$ denote the indicator function defined by
\begin{align*}
s \mapsto \mathbf{1}_{I_l}(s) = \begin{cases} 1 & \text{ if } s \in I_l\\ 0 & \text{ else }. \end{cases}
\end{align*}
Then we take input functions $u$ of the form 
\begin{align*}
u(T-s)= \sum_{l=0}^{N-1} u_l \, \mathbf{1}_{I_l}(T-s)= \sum_{l=0}^{N-1} u_{N-1-l} \, \mathbf{1}_{I_l}(s).
\end{align*}
If $\lambda_k(\theta)\neq 0$ the $k$th component of the solution is then
\begin{align*}
\begin{split}
\tilde \varphi_k(N\tau,u,0)(\theta) &= \int_0^{N\tau} \operatorname{e}^{ \lambda_k(\theta)s}  u(N\tau-s) \di{s} \\
&=   \sum_{l=0}^{N-1} \int_{l\tau}^{(l+1)\tau} \!\!\operatorname{e}^{\lambda_k(\theta)s}  u_{N-1-l} \mathbf{1}_{I_l}(s) \di{s}\\
&=  \left(\frac{ \operatorname{e}^{\tau \lambda_k(\theta)} -1}{\tau \lambda(\theta)}\right)\sum_{l=0}^{N-1} \tau u_{N-1-l} \operatorname{e}^{l \tau \lambda_k(\theta)} .
\end{split}
\end{align*}
If $\lambda_k(\theta)=0$, we have
\begin{align*}
\tilde \varphi_k(N\tau,u,0)=   \sum_{l=0}^{N-1} \int_{l\tau}^{(l+1)\tau} \!\!u_{N-1-l} \mathbf{1}_{I_l}(l\tau + s) \di{s}
=  \sum_{l=0}^{N-1}   \tau \,u_{N-1-l} .
\end{align*}

\medskip
\noindent
{\it Step~3: Approximate the solution by a polynomial}\\
Observing that 
\begin{align*}
\sum_{l=0}^{N-1}  u_{N-1-l} \left(\operatorname{e}^{\tau \lambda_k(\theta)}\right)^l 
\end{align*}
defines a polynomial whose coefficients are given by the values of the input function $u$ defined in \eqref{def:u}, {we define}
\begin{align}\label{eq:def-p}
p(z) := \sum_{l=0}^{N-1} u_{N-1-l} \, z^l.
\end{align}
In terms of {the} polynomial $p${,} the $k$th component of the solution reads as
\begin{align}\label{eq:component_k}
\tilde\varphi_k(N\tau,u,0)=  
\begin{cases}
 \left(\frac{ \operatorname{e}^{\tau \lambda_k(\theta)} -1}{\tau \lambda_k(\theta)}\right) \,\tau \, p\left(\operatorname{e}^{\tau \lambda_k(\theta)}\right) & \quad \text{ if }\lambda_k(\theta)\neq 0\\
 \tau \,  \, p\left(1\right) & \quad \text{ if } \lambda_k(\theta)=0.
\end{cases}
\end{align}
Moreover, since
\begin{align*}
\lim_{\tau \to 0} \frac{\operatorname{e}^{\tau z} -1}{ \tau z} = 1 \quad \text{ for all } z \in \mathbb{C} \setminus\{ 0\},
\end{align*}
there is a $\tau_1>0$ so that for any $\tau \in (0,\tau_1)$ we have 
\begin{align}\label{eq:exp-est_1}
\left| \frac{\operatorname{e}^{\tau \lambda_k(\theta)} -1}{ \tau \lambda_k(\theta)} -1\right|< \tfrac{\varepsilon}{2} \qquad \text{ for all }\, \, \lambda_k(\theta)\neq 0, \, \, k=1,...,n. 
\end{align}
Then, by \eqref{eq:component_k} and \eqref{eq:exp-est_1}{,} for every $\tau \in (0, \tau_1)$ each component of the solution satisfies
\begin{align}\label{eq:exp-est-2}
\left| \tau p(\operatorname{e}^{\tau \lambda_k(\theta)})-\varphi_k(N\tau,u,0)(\theta)  \right|< \tfrac{\varepsilon}{2}.
\end{align}
The significance of \eqref{eq:exp-est-2} is that it {is independent of} the number of input values $u_0,...,u_{N-1}$.

\medskip
\noindent
{\it Step~4: Approximation of each component of the target states by a polynomial}\\
First, we note that there is a $\tau_2>0$ such that the mappings $$\theta \mapsto \operatorname{e}^{\tau\, \lambda_k(\theta)}$$ are injective for all $k=1,...,n$ and all $\tau \in (0,\tau_2)$. Thus, we fix some $\tau <\min\{\tau_1,\tau_2\}$. {The} sets
$$\Omega_k:=\{ \operatorname{e}^{\tau \lambda_k(\theta)}: \theta \in \mathbf{P}\} \subset \mathbb{C}$$
define Jordan arcs. Moreover, let $g(\theta) =T(\theta)f(\theta) $ denote the transformed family of terminal states{, we} consider the continuous function
\begin{align*}
\tilde g_k \colon \Omega_k \to \mathbb{C}, \quad \tilde g_k(z) =  \tfrac{1}{\tau} g_k \left( \lambda^{-1}_k \left( \tfrac{\ln z}{\tau}\right) \right).
\end{align*}
{Applying the results from Section~\ref{sec:approx}, we get} polynomials $p_1,...,p_n$ such that {for each $k=1,...,n$ we have
\begin{align*}
|\tilde g_k(z) - p_k(z)| < \tfrac{\varepsilon}{6\,\tau\, \|T\|_{M,\infty}} \qquad \forall \, \,z \in \Omega_k.
\end{align*}
}
Note that, the latter is equivalent to
\begin{align}\label{eq:est-mergelyan}
| g_k(\theta) - \tau p_k(\operatorname{e}^{\tau \lambda_k(\theta)} )| < \tfrac{\varepsilon}{6\, \|T\|_{M,\infty}} \qquad \forall \, \,\theta \in \p.
\end{align}

\medskip
\noindent
{\it Step~5: Construction of a single polynomial}\\
Let $\Omega = \bigcup_{k=1}^n \Omega_k$. Note that the assumptions imply that the sets $\Omega_1,..., \Omega_n$ are  disjoint. Consider the holomorphic functions $h_1,...,h_n$ defined by
\begin{align*}
h_k\colon \Omega \to \C , \quad   h_k(z)= 
 \begin{cases}
  1 & \text{ if } z \in \Omega_k \\
  0 & \text{ if } z \in \Omega \setminus \Omega_k 
 \end{cases}
\end{align*}
and compute  via Runge's little Theorem polynomials $q_1,...,q_n$ such that
\begin{align*}
\sup_{z \in \Omega} | q_k(z) - h_k(z))| < \frac{\varepsilon}{6\, \tau \, \|T\|_{M,\infty} \, \sum_{l=1}^n \alpha_{k,l}}
\end{align*}
for all $k=1,...,n$, where $\alpha_{k,l}:= \sup_{\theta \in \p}|p_l({\rm e}^{\tau\lambda_k(\theta)})|$.
{Then} we define the polynomial
\begin{equation}\label{eq:def_poly_cont_si_case}
p(z) = \sum_{k=1}^n p_k(z)q_k(z).
\end{equation}

\medskip
\noindent
{\it Step~6: Definition of the input values}\\
Determine a monomial representation of the polynomial $p$ in \eqref{eq:def_poly_cont_si_case}, i.e.{,}
\begin{equation}\label{eq:def_poly_cont_si_case_monic}
p(z) = \sum_{k=1}^n p_k(z)q_k(z) = \sum_{k=0}^K  c_k z^k
\end{equation}
and set $N=K+1$ and define the input $u \colon [0,N\tau]\to \C$ by 
\begin{align*}
 u|_{I_l} (t) := c_l, \quad l=0,...,N-1.
\end{align*}
Then, putting things together, for given $f \in C_n(\p)$ and $\varepsilon>0$  we get for each component 
\begin{align*}
|\varphi_k(N\tau,u,0)(\theta)- f_k(\theta)|  &\leq \|T\|_{M,\infty} \, \left|\tilde \varphi_k(N\tau,u,0)(\theta)-\tau p(\operatorname{e}^{\tau \lambda_k(\theta)} ) \right| 
\\ & \qquad \qquad + \|T\|_{M,\infty} \, \left| \tau p(\operatorname{e}^{\tau \lambda_k(\theta)} )- g_k(\theta) \right| \\
&< \tfrac{\e}{2} +  \|T\|_{M,\infty}\,  \left|  g_k(\theta) - \tau p_k(\operatorname{e}^{\tau \lambda_k(\theta)} ) q_k( \operatorname{e}^{\tau \lambda_k(\theta)}) \right| \\
& \qquad \qquad+ \|T\|_{M,\infty} \,\tau \,\sum_{l\neq k}  \left| p_l(\operatorname{e}^{\tau \lambda_k(\theta)}) \right| \,\left|q_l(\operatorname{e}^{\tau \lambda_k(\theta)}) \right|\\
&< \tfrac{\e}{2} + \tfrac{\e}{6} +  \|T\|_{M,\infty} \,  \left|  g_k(\theta) - \tau p_k(\operatorname{e}^{\tau \lambda_k(\theta)} ) \right| \\
& \qquad \qquad+ \|T\|_{M,\infty} \, \tau \left|  1-q_k( \operatorname{e}^{\tau \lambda_k(\theta)})\right|\, \left| p_k(\operatorname{e}^{\tau \lambda_k(\theta)} ) \right|\\
& < \tfrac{\e}{2} + \tfrac{\e}{6} + \tfrac{\e}{6} + \tfrac{\e}{6} = \e .
\end{align*}
Consequently, we have
\begin{align*}
\|f - \varphi(N\tau,u,0)\|_{\infty} < \varepsilon.
\end{align*}

\begin{remark}
The approach in Section~\ref{sec:poly_approx_cont} was used in a first and direct proof of the sufficiency of (N1), (N2) and (S2) for uniform ensemble reachability for single-input systems, cf. \cite{helmke2014uniform}. Indeed, since the set $\Omega$ is compact with empty interior and $\C \setminus \Omega$ is connected, it is possible to applying Mergelyan's theorem directly to the continuous function
\begin{align*}
g \colon \Omega \to \mathbb{C}, \quad g{|_{\Omega_k}}(z) =  \tfrac{1}{\tau} g_k \left( \lambda^{-1}_k \left( \tfrac{\ln z}{\tau}\right) \right).
\end{align*}
Thus, there is a polynomial $p$ such that
\begin{align*}
|g(z) - \tau p(z)|< \e \qquad \forall \, z \in \Omega.
\end{align*}
Also we note that, the property that the interior of $\Omega$ is empty is a special case of Mergelyan's result, which was proven earlier by Lavrientev in 1934, cf. \cite{gaier1987}. However, Lavrientev's proof is not available to the author. Maybe the proof can be used to obtain another constructive method{, at least for} special cases.
\end{remark}

\section{Comments}\label{sec:comments}

In this section we discuss some perspectives on the methods presented in Section~\ref{sec:poly_si_methods}.

\medskip
\noindent
{\it Degrees of freedom}

If the approximation domains are given by Jordan arcs, the methods for (S1) and (S2) have some degrees of freedom. For instance, there are plenty of ways to close a Jordan arc to a Jordan curve. This in turn frees up space for designing the conformal mappings such that the conditions \eqref{eq:conf-cond_p-1}, \eqref{eq:conf-cond_p-2} and \eqref{eq:conf-cond_p,k-1}, \eqref{eq:conf-cond_p,k-2} are satisfied. Thus, the computation of the conformal mappings and the design of the Jordan curves can be treated simultaneously. Besides, for pairs $(A,b)$ where $A$ is block diagonal and {each} block satisfies {the} sufficiency condition (S1) or (S2), the methods can be applied for every block and an overall polynomial can be constructed {in the same manner} the polynomials $q_1,...,q_n$ are determined in the  methods for (S2).

\medskip
\noindent
{\it Real approximation domains}

In case the domains, where the approximations take place{,} are compact real intervals, the construction process is much easier. This can be {achieved} by considering a mixture of open-loop and feedback control inputs of the form
\[ K(\theta)x_t(\theta) + u_t. \]
Indeed, due to the necessary condition (N1), {a suitable feedback term $K(\theta)x_t(\theta)$ can be used to shift the spectra of $A(\theta)$ to an appropriate position}. That is, via $K(\theta)$ the matrices $A(\theta)+B(\theta)K(\theta)$ can have any spectra. Equivalently, the coefficients of the characteristic polynomials of $A(\theta)+B(\theta)K(\theta)$ can be arbitrarily designed. For more details on this we refer to \cite[Section~6.4]{fuhrmann2015mathematics}, \cite{laa_feedback}.

\medskip
\noindent
{\it Multi-input systems}

The known necessary and sufficient conditions for ensemble reachability are in the multi-input case much less precise compared to the single-input systems. For special classes of systems, such as $(A(\theta),B(\theta))$ upper triangular, the problem can be traced back to the single-input case. In this context, a general but very strong sufficiency condition is that the Hermite indices of the pair $(A(\theta),B(\theta))$ are constant. Technicalities aside, these systems can be transformed into the Hermite canonical form 
\begin{equation*}
\begin{pmatrix}
\tilde A_{11}(\theta)  & \cdots & \tilde A_{1k}(\theta) \\
 & \ddots&  \vdots\\
0 & & \tilde A_{kk}(\theta)
\end{pmatrix},
\quad
\begin{pmatrix}
\tilde b_{1}(\theta) &  & 0 & * & * \\
 & \ddots & &*& * \\
0 & & \tilde b_{k}(\theta)&*& * 
 \end{pmatrix}, 
\end{equation*}
and the construction problem can be solved by applying the single-input  methods to the subpairs $(\tilde A_{ii}(\theta),\tilde b_i(\theta))$. We note that it is not required to use one of the methods for all subpairs. In particular, this applies to the often studied class $(\theta A,B)$. For more details{,} we refer to \cite{JDE_ensembles_2021} and the references therein.

\medskip
\noindent
{\bf Acknowledgements}
The author is grateful to Vladimir Andrievskii and Oliver Roth  for many comments and valuable suggestions on complex approximation and conformal mappings. This research was supported by the German Research Foundation (DFG) within the grant SCHO 1780/1-1 and the European Research Council (ERC) project CHRiSHarMa no. DLV-682402.

%

\end{document}